\documentclass[10pt]{amsart}
\usepackage{pgf,tikz}
\usetikzlibrary{decorations.pathreplacing}

\usepackage{caption}

\usepackage{amssymb}
\usepackage{stmaryrd}
\usepackage{mathrsfs}
\usepackage[all]{xy}
\makeatletter
  \def\@wrindex#1{%
    \protected@write\@indexfile{}%
      {\string\indexentry{#1}{ \S\thesubsection (p.\thepage)}}
    \endgroup
  \@esphack}

\makeatother

\makeindex

\newcommand{\BB}{\rho}
\newcommand{\DR}{\mathrm{DR}}
\newcommand{\B}{\mathrm{B}}
\newcommand{\KZ}{\mathrm{KZ}}

\usepackage{rotating}
\usepackage{amscd}
\usepackage{epsfig}

\author{Benjamin Enriquez}

\address{Institut de Recherche Math\'{e}matique Avanc\'{e}e, UMR 7501, 
Universit\'{e} de Strasbourg et CNRS, 7
rue Ren\'{e} Descartes, 67000 Strasbourg, France}
\email{enriquez@math.unistra.fr}
\date{May 12th, 2020.}

\newtheorem{thm}{Theorem}[section]
\newtheorem{lem}[thm]{Lemma}

\newtheorem{prop}[thm]{Proposition}

{\theoremstyle{definition} \newtheorem{rem}[thm]{Remark}}
{\theoremstyle{definition} \newtheorem{defn}[thm]{Definition}}

{\theoremstyle{remark} }

\numberwithin{equation}{subsection}
\numberwithin{figure}{section}

\begin{document}

\baselineskip 16pt 

\title[The Betti side of the double shuffle theory: a survey]{The Betti side of the double shuffle theory: a survey}

\begin{abstract}
This is a survey of \cite{EF1,EF2,EF3}. The purpose of this series of papers is: (1) to give a proof that associator relations imply 
double shuffle relations, alternative to \cite{F3}; (2) to make explicit the bitorsor structure on Racinet's torsor of double shuffle 
relations. The main tool is the interpretation of the harmonic coproduct in terms of the topology of the moduli spaces 
$\mathfrak M_{0,4}$ and $\mathfrak M_{0,5}$, introduced in \cite{DeT}, and its extension to the Betti setup. 
\end{abstract}

\bibliographystyle{amsalpha+}
\maketitle


{\footnotesize \tableofcontents}

\section*{Introduction}

The multizeta values (MZVs) are the real numbers defined by the series
$$
\zeta(k_1,\cdots,k_m):=\sum_{n_1>\cdots>n_m>0}\frac{1}
{n_1^{k_1}\cdots n_m^{k_m}}
$$
for $k_1$,\dots, $k_m\in{\mathbb Z}_{>0}$ and $k_1>1$. 
These numbers have recently garnered much interest due to their appearance in various fields of physics and mathematics
(\cite{BK}). They appear to be examples of periods (\cite{KZ}) and are, as such, related with motive theory (\cite{De}). 
Using this theory, upper bounds for dimensions of spaces of MZVs have been obtained (\cite{DG,Te}). A related problem is the identification
of the algebraic and linear relations among MZVs.

A recent review of the available systems of relations can be found in \cite{Bach}. Among them, we will focus of the interrelations 
between: (a) the associator system of relations (\cite{Dr,LM}); (b) the regularized double shuffle relations (\cite{IKZ, Rac}).   

Each system of relations gives rise to a $\mathbb Q$-scheme, defined as the spectrum of the free commutative 
$\mathbb Q$-algebra over formal variables $\zeta^f(k_1,\cdots,k_m)$ for $(k_1,\ldots, k_m)\in{\mathbb Z}_{>1}
\times({\mathbb Z}_{>0})^{m-1}$ by the corresponding ideal. These schemes are called the {\it scheme of associators} 
in case (a) and the {\it double shuffle scheme} in case (b); for $\mathbf k$ a $\mathbb Q$-algebra, the sets of 
$\mathbf k$-points of these schemes are denoted $M(\mathbf k)$ in case (a) and $\mathsf{DMR}^{\DR,\B}(\mathbf k)$
in case (b). The definition of $M(\mathbf k)$ can be found in \cite{Dr}, p. 848, and $\mathsf{DMR}^{\DR,\B}(\mathbf k)=\sqcup_{\mu\in\mathbf k^\times}\mathsf{DMR}_\mu(\mathbf k)$, where $\mathsf{DMR}_\mu(\mathbf k)$ is as in \cite{Rac}, 
Déf. 3.2.1. The structures of these schemes is elucidated by the following results. 

\begin{thm} (\cite{Dr})
(1) There are explicit $\mathbb Q$-group schemes $\{\mathbb Q$-algebras$\}\ni \mathbf k\mapsto \mathrm{GT}(\mathbf k)$, 
$\mathrm{GRT}(\mathbf k)\in\{$groups$\}$, and for any $\mathbb Q$-algebra $\mathbf k$, commuting left and right free and 
transitive actions of $\mathrm{GT}(\mathbf k)$ and $\mathrm{GRT}(\mathbf k)$ on $M(\mathbf k)$. 
 
 (2) These $\mathbb Q$-group schemes are extensions of the multiplicative group $\mathbb G_m$ by prounipotent group schemes. 
Their Lie algebras $\mathfrak{gt}$ and $\mathfrak{grt}$ are filtered, moreover $\mathfrak{grt}$ is complete graded.  
\end{thm} 

Note that the group $\mathrm{Aut}_{\mathrm{GRT}(\mathbf k)}(M(\mathbf k))$ of permutations of $M(\mathbf k)$ which 
commute with the action of $\mathrm{GRT}(\mathbf k)$ naturally acts on $M(\mathbf k)$. (1) says that there is an isomorphism 
between this group and the explicit group $\mathrm{GT}(\mathbf k)$, which is compatible with their actions on $M(\mathbf k)$. 

\begin{thm} (\cite{Rac}) (1) There is an explicit  $\mathbb Q$-group scheme $\{\mathbb Q$-algebras$\}\ni \mathbf k\mapsto 
\mathsf{DMR}^\DR(\mathbf k)\in\{$groups$\}$, and for any $\mathbb Q$-algebra $\mathbf k$, a free and transitive left action 
of $\mathsf{DMR}^\DR(\mathbf k)$ on $\mathsf{DMR}^{\DR,\B}(\mathbf k)$. 

(2) This $\mathbb Q$-group scheme is an extension of the multiplicative group $\mathbb G_m$ by  a prounipotent group scheme. 
Its Lie algebra $\mathfrak{dmr}^\DR$ is complete graded.  
\end{thm} 

This formulation is obtained in \cite{EF2} using the main result of \cite{Rac}. The best available result on the comparison of 
the associator and double shuffle schemes is as follows. 

\begin{thm} \label{thm:0:3} (\cite{F3}) 
For any $\mathbb Q$-algebra $\mathbf k$, there are compatible inclusions of sets $M(\mathbf k)\subset 
\mathsf{DMR}^{\DR,\B}(\mathbf k)$ and of groups $ \mathrm{GRT}(\mathbf k)^{\mathrm{op}}\subset
\mathsf{DMR}^\DR(\mathbf k)$ (where $\mathrm{op}$ denotes the opposite group).  
\end{thm} 

The proof in \cite{F3} relies on the construction, out of the family of multiple polylogarithm functions, of elements in the bar-complex 
of the moduli space $\mathfrak M_{0,5}$, which are then viewed as linear forms on the enveloping algebra $U(\mathfrak p_5)$
(see \S\ref{sect:2:3:27032020}), and on the study of the combinatorics of these linear forms. This result was also announced in 
the unfinished preprint \cite{DeT}, which contains in particular a description of one of the main actors of double shuffle theory, 
the `harmonic coproduct', in terms of topology of the moduli spaces $\mathfrak M_{0,4}$ and $\mathfrak M_{0,5}$. 

 The main objectives of the series of papers \cite{EF1,EF2,EF3} are: (a) giving a new proof of Theorem \ref{thm:0:3}, based on the ideas of \cite{DeT} (\cite{EF2}); (b) making the group $\mathrm{Aut}_{\mathsf{DMR}^{\DR}(\mathbf k)}(\mathsf{DMR}^{\DR,\B}
(\mathbf k))$ explicit, together with its action on $\mathsf{DMR}^{\DR,\B}(\mathbf k)$ (\cite{EF3}). In order to reach them, 
we perform an intermediate task: (c) constructing a `Betti' version of the algebraic apparatus of double shuffle theory and showing 
how it is related to the original one by an associator (\cite{EF1}). The material relative to objective (c) (resp., (a), (b)) is reviewed in 
\S\ref{sect:1:27032020} (resp. \S\ref{sect:2:27032020},\S\ref{sect:3:27032020}).

\section{The algebraic framework of the double shuffle theory}\label{sect:1:27032020}

\subsection{The de Rham side of double shuffle theory} 

Let $\mathcal V^\DR$ be the free associative $\mathbf k$-algebra over generators $e_0,e_1$; it is 
$\mathbb Z_{\geq0}$-graded, with $e_0,e_1$ being of degree 1. Let 
$\mathcal W^\DR:=\mathbf k1\oplus\mathcal V^\DR e_1$; this is a $\mathbb Z_{\geq0}$-graded
subalgebra of $\mathcal V^\DR$. Set $\mathcal M^\DR:=\mathcal V^\DR/\mathcal V^\DR e_0$; this 
is a $\mathbb Z_{\geq0}$-graded left $\mathcal V^\DR$-module, therefore by restriction a left 
$\mathcal W^\DR$-module, which is free of rank one, generated by the class $1_\DR\in\mathcal M^\DR$ 
of the element $1\in\mathcal V^\DR$.

Let $\Delta^{\mathcal V,\DR}:\mathcal V^\DR\to(\mathcal V^\DR)^{\otimes2}$ be the $\mathbf k$-algebra morphism 
such that $e_i\mapsto e_i\otimes1+1\otimes e_i$ for $i=0,1$. One shows that $\mathcal W^\DR$ is freely generated, 
as an associative algebra, by its elements $y_n:=-e_0^{n-1}e_1$, where $n\geq1$. We denote by 
$\Delta^{\mathcal W,\DR}:\mathcal W^\DR\to(\mathcal W^\DR)^{\otimes2}$ the $\mathbf k$-algebra morphism 
such that $y_n\mapsto\sum_{i=0}^n y_i\otimes y_{n-i}$ for $n\geq1$, where $y_0:=1$, and by 
$\Delta^{\mathcal M,\DR}:\mathcal M^\DR\to(\mathcal M^\DR)^{\otimes2}$ the $\mathbf k$-module morphism 
such that $\Delta^{\mathcal M,\DR}(a\cdot 1_\DR)=\Delta^{\mathcal W,\DR}(a)\cdot 1_\DR^{\otimes2}$, where 
$\cdot$ denotes the action of $\mathcal V^\DR$ on $\mathcal M^\DR$. The maps $\Delta^{\mathcal X,\DR}$, 
$\mathcal X\in\{\mathcal V,\mathcal W,\mathcal M\}$ are all compatible with the $\mathbb Z_{\geq0}$-gradings. 

Then $(\mathcal V^\DR,\Delta^{\mathcal V,\DR})$ and $(\mathcal W^\DR,\Delta^{\mathcal W,\DR})$ are cocommutative 
Hopf algebras, but the inclusion $\mathcal W^\DR\subset\mathcal V^\DR$ is not compatible with the coproducts; 
$(\mathcal M^\DR,\Delta^{\mathcal M,\DR})$  is a cocommutative coalgebra, and is a coalgebra-module over 
$(\mathcal W^\DR,\Delta^{\mathcal W,\DR})$. 

We will denote by $\hat X$ (or $X^\wedge$) the completion of a $\mathbb Z_{\geq0}$-graded $\mathbf k$-module $X$, 
and use the same notation for the completion of a morphism between such objects. 

For $g\in\hat{\mathcal V}^\DR$, let $\Gamma_g(t):=\mathrm{exp}(\sum_{n\geq1}(-1)^{n+1}(g|e_0^{n-1}e_1)t^n/n)
\in\mathbf k[[t]]^\times$, where $w\mapsto (g|w)$ is the map $\{e_0,e_1\}^*\to\mathbf k$ such that $g=\sum_{w\in
\{e_0,e_1\}^*}(g|w)w$. 

Let $w\mapsto w_{\mathrm{reg}}$ be the composed map $e_0\mathcal V^\DR e_1\hookrightarrow\mathcal V^\DR
\to\mathcal V^\DR\otimes_{\mathbb Q}\mathbb Q[\alpha_0,\alpha_1]\to\mathcal V^\DR$, where the second map is the 
$\mathbf k$-algebra morphism induced by $e_i\mapsto e_i\otimes 1-1\otimes\alpha_i$ for $i=0,1$ and the third map is 
the $\mathbf k$-module map induced by $v\otimes \alpha_0^a\alpha_1^b\mapsto e_1^b v e_0^a$ for $a,b\geq 0$, 
$v\in\mathcal V^\DR$. For $w\in e_0\{e_0,e_1\}^*e_1$, let $m(w)$ be the number of occurences of $e_1$ in $w$.  

Let $\zeta:e_0\{e_0,e_1\}^*e_1\to\mathbb R$ be the map defined by $\zeta(e_0^{k_1}e_1\cdots e_0^{k_m}e_1)
=\zeta(k_1,\ldots,k_m)$ for $m\geq 1$, $(k_1,\ldots,k_m)\in\mathbb Z_{>1}\times(\mathbb Z_{>0})^{m-1}$. 
Set\footnote{The notation $\hat{\mathcal V}_{\mathbb C}^\DR$, $\hat{\mathcal M}_{\mathbb C}^\DR$ stands for 
the specializations of $\hat{\mathcal V}^\DR$, $\hat{\mathcal M}^\DR$ for $\mathbf k=\mathbb C$.} 
$$
\varphi_{\mathrm{KZ}}:=\sum_{w\in e_0\{e_0,e_1\}^*e_1}(-1)^{m(w)}\zeta(w)w_{\mathrm{reg}}\in
\hat{\mathcal V}_{\mathbb C}^\DR;  
$$
this is a generating series for the MZVs, usually called the {\it Knizhnik-Zamolodchikov associator} (see \cite{Dr,LM,F1}).

The system of regularization and double shuffle relations (\cite{IKZ,Rac}) between MZVs can be formulated as follows:  
\begin{equation}\label{DS:PhiKZ:1}
\varphi_{\mathrm{KZ}}\in\mathcal G(\hat{\mathcal V}_{\mathbb C}^\DR),\quad 
(\Gamma_{\varphi_{\mathrm{KZ}}}(-e_1)^{-1}\varphi_{\mathrm{KZ}})\cdot 1_\DR\in
\mathcal G(\hat{\mathcal M}_{\mathbb C}^\DR),
\end{equation}
\begin{equation}\label{DS:PhiKZ:2}
(\varphi_{\mathrm{KZ}}|e_0)=(\varphi_{\mathrm{KZ}}|e_1)=0, 
(\varphi_{\mathrm{KZ}}|e_0e_1)=(2\pi i)^2/24, 
\end{equation}
where $\mathcal G$ denotes the set of group-like elements of $\hat{\mathcal V}^\DR$ (resp. $\hat{\mathcal M}^\DR$) 
for $\hat\Delta^{\mathcal V,\DR}$  (resp. $\hat\Delta^{\mathcal M,\DR}$). 

This formulation leads to the following definition: 
\begin{defn}
For $\mu\in\mathbf k$, one defines $\mathsf{DMR}_\mu(\mathbf k)$ to be the set of elements $\Phi\in\hat{\mathcal V}^\DR$
satisfying relations \eqref{DS:PhiKZ:1}, \eqref{DS:PhiKZ:2}, with $\varphi_{\mathrm{KZ}}$, $2\pi i$, $\mathbb C$
replaced by $\Phi$, $\mu$, $\mathbf k$. 
\end{defn}
One has therefore $\varphi_{\mathrm{KZ}}\in\mathsf{DMR}_{2\pi i}(\mathbb C)$.  

\begin{rem} The `double shuffle' system of relations is the conjunction of the systems of harmonic and shuffle relations, 
which follow from the expressions of the MZVs respectively as iterated sums and as iterated integrals. 
An example of a harmonic relation is 
$$
\forall a,b>1,\quad \zeta(a)\zeta(b)=\zeta(a+b)+\zeta(a,b)+\zeta(b,a), 
$$
which follows from $\zeta(a)\zeta(b)=\sum_{n,m>0}n^{-a}m^{-b}
=(\sum_{n=m>0}+\sum_{n>m>0}+\sum_{m>n>0})n^{-a}m^{-b}
=\zeta(a,b)+\zeta(a,b)+\zeta(b,a)$. 
An example of a shuffle relation is 
$$
\forall a,b>1,\quad \zeta(a)\zeta(b)=\sum_{i+j=a+b}\Big(
\begin{pmatrix}a-1 \\ i-1\end{pmatrix}
+\begin{pmatrix}b-1 \\ j-1\end{pmatrix}\Big)\zeta(i,j), 
$$
which follows from $\zeta(c)=\int_{0<s_1<\ldots< s_c<1}{{ds_1}\over{1-s_1}}\wedge{{ds_2}\over{s_2}}
\wedge\cdots\wedge{{ds_c}\over{s_c}}$, $\zeta(i,j)=\int_{0<s_1<\ldots< s_{i+j}<1}
{{ds_1}\over{1-s_1}}\wedge{{ds_2}\over{s_2}}\wedge\cdots\wedge{{ds_j}\over{s_j}}\wedge 
{{ds_{j+1}}\over{1-s_{j+1}}}\wedge{{ds_{j+2}}\over{s_{j+2}}}\wedge\cdots\wedge{{ds_{i+j}}\over{s_{i+j}}}$ 
and the shuffle identity for products of iterated integrals. 
\end{rem}

\begin{rem}
The relation between the formalism of \cite{EF1} and \cite{Rac} is as follows: the elements $e_0,e_1$
correspond to $x_0,-x_1$; the pair $(\hat{\mathcal V}^\DR,\hat\Delta^{\mathcal V,\DR})$ corresponds to 
$(\mathbf k\langle\langle X\rangle\rangle,\hat\Delta)$; the pairs $(\hat{\mathcal W}^\DR,\hat\Delta^{\mathcal W,\DR})$ 
and $(\hat{\mathcal M}^\DR,\hat\Delta^{\mathcal M,\DR})$ both correspond to 
$(\mathbf k\langle\langle Y\rangle\rangle,\hat\Delta_\star)$; the map $\hat{\mathcal V}^\DR\to\hat{\mathcal M}^\DR$, 
$a\mapsto a\cdot 1_\DR$ corresponds to $\pi_Y$. 
\end{rem}

\subsection{The Betti side of double shuffle theory}

Let $\mathcal V^\B$ be the $\mathbf k$-algebra with generators $X_0^{\pm1},X_1^{\pm1}$ and relations 
$X_iX_i^{-1}=X_i^{-1}X_i=1$ for $i=0,1$. It is equipped with the filtration $\mathcal V^\B=F^0\mathcal V^\B
\supset F^1\mathcal V^\B\supset\cdots$, where $F^k\mathcal V^\B:=(\mathcal V^\B_+)^k$ for $k\geq0$, where 
$\mathcal V^\B_+$ is the kernel of the $\mathbf k$-algebra morphism $\mathcal V^\B\to\mathbf k$ given by 
$X_i^{\pm1}\mapsto 1$ for $i=0,1$. 

Define a $\mathbf k$-subalgebra $\mathcal W^\B$ of $\mathcal V^\B$ by $\mathcal W^\B:=\mathbf k1\oplus
\mathcal V^\B(X_1-1)$. It is equipped with the induced filtration $F^k\mathcal W^\B:=\mathcal W^\B\cap F^k\mathcal V^\B$
for $k\geq0$. One can show that $\mathcal W^\B$ is presented by generators $X_1^{\pm1}$ and 
$Y_n^{\pm}:=(X_0^{\pm1}-1)^{n-1}X_0^{\pm1}(1-X_1^{\pm1})$ for $n>0$, and relations
$X_1X_1^{-1}=X_1^{-1}X_1=1$ (see \cite{EF1}, \S2.2). 

Set $\mathcal M^\B:=\mathcal V^\B/\mathcal V^\B(X_0-1)$. This is a left $\mathcal V^\B$-module, hence by restriction 
a left $\mathcal W^\B$-module; we denote by $1_\B\in\mathcal M^\B$ the class of $1\in\mathcal V^\B$ and by $\cdot$
the action of $\mathcal V^\B$ on $\mathcal M^\B$. Then $\mathcal M^\B$ is free of rank one as a $\mathcal W^\B$-module, 
generated by $1_\B$. One equips $\mathcal M^\B$ with the filtration $F^k\mathcal M^\B:=F^k\mathcal V^\B\cdot 1_\B$ 
for $k\geq0$. 

The filtrations of $\mathcal V^\B$, $\mathcal W^\B$ and $\mathcal M^\B$ are compatible with the algebra inclusion 
$\mathcal V^\B\subset\mathcal W^\B$ and the algebra actions of $\mathcal V^\B$ and $\mathcal W^\B$ on $\mathcal M^\B$. 

There is a unique $\mathbf k$-algebra morphism $\Delta^{\mathcal V,\B}:\mathcal V^\B\to(\mathcal V^\B)^{\otimes2}$, such that 
$X_i^{\pm1}\mapsto (X_i^{\pm1})^{\otimes2}$ for $i=0,1$. There is a unique $\mathbf k$-algebra morphism 
$\Delta^{\mathcal W,\B}:\mathcal W^\B\to(\mathcal W^\B)^{\otimes2}$, such that $X_1^{\pm1}\mapsto 
(X_1^{\pm1})^{\otimes2}$ and $Y_k^\pm\mapsto\sum_{i=0}^k Y_i^\pm\otimes Y_{k-i}^\pm$, where
$Y_0^\pm=1$. There is a unique $\mathbf k$-module morphism 
$\Delta^{\mathcal M,\B}:\mathcal M^\B\to(\mathcal M^\B)^{\otimes2}$, such that 
$\Delta^{\mathcal M,\B}(a\cdot 1_\B)=\Delta^{\mathcal W,\B}(a)\cdot 1_\B^{\otimes2}$ for any $a\in\mathcal W^\B$. 

As before, $(\mathcal V^\B,\Delta^{\mathcal V,\B})$ and $(\mathcal W^\B,\Delta^{\mathcal W,\B})$ are cocommutative 
Hopf algebras, and the inclusion $\mathcal W^\B\subset\mathcal V^\B$ is not compatible with the coproducts; 
$(\mathcal M^\B,\Delta^{\mathcal M,\B})$  is a cocommutative coalgebra, and is a coalgebra-module over 
$(\mathcal W^\B,\Delta^{\mathcal W,\B})$. 

The coproducts $\Delta^{\mathcal X,\B}$, with $\mathcal X\in\{\mathcal V,\mathcal W,\mathcal M\}$ are all compatible with 
the filtrations. We denote by $\hat X$ (or $X^\wedge$) the completion of a filtered $\mathbf k$-module, and use the same 
notation for the completion of a morphism between such objects. 

\subsection{Filtrations and gradings} 

As $(\mathcal V^\B,\Delta^{\mathcal V,\B})$, $(\mathcal W^\B,\Delta^{\mathcal W,\B})$ and 
$(\mathcal M^\B,\Delta^{\mathcal M,\B})$ are Hopf algebras and a coalgebra in the category of filtered
$\mathbf k$-modules, the associated graded objects have the same status in the category of 
$\mathbb Z_{\geq0}$-modules; these objects are  respectively isomorphic to 
$(\mathcal V^\DR,\Delta^{\mathcal V,\DR})$, $(\mathcal W^\DR,\Delta^{\mathcal W,\DR})$ and 
$(\mathcal M^\DR,\Delta^{\mathcal M,\DR})$. The isomorphism $\mathrm{gr}(\mathcal V^\B)\simeq\mathcal V^\DR$
is induced by $\mathrm{gr}_1(\mathcal V^\B)\ni(\text{class of }X_i-1)\mapsto e_i\in\mathcal V^\DR$ for $i=0,1$; 
it induces the isomorphism $\mathrm{gr}(\mathcal W^\B)\simeq\mathcal W^\DR$. The isomorphism 
$\mathrm{gr}(\mathcal M^\B)\simeq\mathcal M^\DR$ is based on the fact that the isomorphism $\mathcal W^\B\to
\mathcal M^\B$, $a\mapsto a\cdot 1_\B$, induces an isomorphism $F^i\mathcal W^\B\to F^i\mathcal M^\B$ for any $i\geq0$, 
which induces the first map in the sequence of isomorphisms $\mathrm{gr}_i(\mathcal M^\B)\simeq\mathrm{gr}_i(\mathcal W^\B)
\simeq\mathcal W^\DR\to\mathcal M^\DR$, the last map being $a\mapsto a\cdot 1_\DR$. 
 
\subsection{Comparison isomorphisms and geometric interpretations}

\subsubsection{Automorphisms of the de Rham side}\label{sect:1:4:1:0305}

Set $G^\DR(\mathbf k):=\mathbf k^\times\times\mathcal G(\hat{\mathcal V}^\DR)$. For $(\mu,g)\in G^\DR(\mathbf k)$, let 
$\mathrm{aut}_{(\mu,g)}^{\mathcal V,\DR,(1)}$ be the topological $\mathbf k$-algebra automorphism of 
$\hat{\mathcal V}^\DR$ defined by $e_0\mapsto  g\cdot\mu e_0\cdot g^{-1}$, $e_1\mapsto \mu e_1$. Let 
$\mathrm{aut}_{(\mu,g)}^{\mathcal V,\DR,(10)}$ be the topological $\mathbf k$-module automorphism 
of $\hat{\mathcal V}^\DR$ defined by $\mathrm{aut}_{(\mu,g)}^{\mathcal V,\DR,(10)}(a):=
\mathrm{aut}_{(\mu,g)}^{\mathcal V,\DR,(1)}(a)\cdot g$ for any $a\in\hat{\mathcal V}^\DR$. One checks that 
$\mathrm{aut}_{(\mu,g)}^{\mathcal V,\DR,(1)}$ restricts to a topological $\mathbf k$-algebra automorphism of 
$\hat{\mathcal W}^\DR$, denoted $\mathrm{aut}_{(\mu,g)}^{\mathcal W,\DR,(1)}$, and that there is a unique 
topological $\mathbf k$-module automorphism $\mathrm{aut}_{(\mu,g)}^{\mathcal M,\DR,(10)}$ of 
$\hat{\mathcal M}^\DR$, such that $\mathrm{aut}_{(\mu,g)}^{\mathcal M,\DR,(10)}(a\cdot 1_\DR)
=\mathrm{aut}_{(\mu,g)}^{\mathcal V,\DR,(10)}(a)\cdot 1_\DR$ for any $a\in\hat{\mathcal V}^\DR$. 

One checks that $(\mu,g)\circledast(\mu',g'):=(\mu\mu',\mathrm{aut}_{(\mu,g)}^{\mathcal V,\DR,(10)}(g'))$
equips $G^\DR(\mathbf k)$ with a group structure, of which $\mathcal G(\hat{\mathcal V}^\DR)$ is a subgroup. 
For $A$ an algebra and $M$ an $A$-module, denote by 
$\mathrm{Aut}(A,M)$ the set of pairs $(\alpha,\theta)$, where $\alpha$ is an algebra automorphism of $A$, and 
$\theta$ is an automorphism of $M$, such that $\theta(am)=\alpha(a)\theta(m)$ for $a\in A$, $m\in M$; this is naturally 
a group. Then the map taking $(\mu,g)$ to 
$(\mathrm{aut}_{(\mu,g)}^{\mathcal W,\DR,(1)},\mathrm{aut}_{(\mu,g)}^{\mathcal M,\DR,(10)})$ is a 
group morphism from $(G^\DR(\mathbf k),\circledast)$ to $\mathrm{Aut}(\hat{\mathcal W}^\DR,\hat{\mathcal M}^\DR)$. 

The map $\Gamma:G^\DR(\mathbf k)\to(\hat{\mathcal W}^\DR)^\times$, $(\mu,g)\mapsto \Gamma_g^{-1}(-e_1)$
satisfies the cocycle identity $\Gamma((\mu,g)\circledast(\mu',g'))=\Gamma(\mu,g)
\mathrm{aut}_{(\mu,g)}^{\mathcal W,\DR,(1)}(\Gamma(\mu',g'))$. It follows that the map taking $(\mu,g)$ to
$({}^\Gamma\!\mathrm{aut}_{(\mu,g)}^{\mathcal W,\DR,(1)},
{}^\Gamma\!\mathrm{aut}_{(\mu,g)}^{\mathcal M,\DR,(10)})$ is a group morphism
from $(G^\DR(\mathbf k),\circledast)$ to $\mathrm{Aut}(\hat{\mathcal W}^\DR,\hat{\mathcal M}^\DR)$, 
where\footnote{If $A$ is an algebra and $u\in A^\times$, then $\mathrm{Ad}_u$ is the automorphism 
of $A$ given by $a\mapsto uau^{-1}$; if moreover $M$ is a left $A$-module, then $\ell_u$ is the automorphism of 
$M$ given by $m\mapsto um$.} 
${}^\Gamma\!\mathrm{aut}_{(\mu,g)}^{\mathcal W,\DR,(1)}:=\mathrm{Ad}_{\Gamma(\mu,g)}\circ
\mathrm{aut}_{(\mu,g)}^{\mathcal W,\DR,(1)}$ and 
${}^\Gamma\!\mathrm{aut}_{(\mu,g)}^{\mathcal M,\DR,(10)}:=
\ell_{\Gamma(\mu,g)}\circ\mathrm{aut}_{(\mu,g)}^{\mathcal M,\DR,(10)}$.

\begin{rem}
For $g\in\mathcal G(\hat{\mathcal V}^\DR)$, the automorphisms $S_g$, $S_{\Theta(g)}$, $S_g^Y$, $S^Y_{\Theta(g)}$
from \cite{Rac} correspond to $\mathrm{aut}_{(1,g)}^{\mathcal V,\DR,(10)}$, 
$^\Gamma\!\mathrm{aut}_{(1,g)}^{\mathcal V,\DR,(10)}$, $\mathrm{aut}_{(1,g)}^{\mathcal M,\DR,(10)}$, 
$^\Gamma\!\mathrm{aut}_{(1,g)}^{\hat{\mathcal M}^\DR,(10)}$.  
 \end{rem}

\subsubsection{Isomorphisms between the Betti and de Rham sides}

There is a topological $\mathbf k$-algebra isomorphism $i^{\mathcal V}:\hat{\mathcal V}^\B\to\hat{\mathcal V}^\DR$
defined by $X_i\mapsto\mathrm{exp}(e_i)$ for $i=0,1$. It restricts to a topological $\mathbf k$-algebra isomorphism 
$i^{\mathcal W}:\hat{\mathcal W}^\B\to\hat{\mathcal W}^\DR$ and it induces a $\mathbf k$-module isomorphism 
$i^{\mathcal M}:\hat{\mathcal M}^\B\to\hat{\mathcal M}^\DR$, defined by $i^{\mathcal M}(a\cdot 1_\DR)
=i^{\mathcal V}(a)\cdot 1_\DR$ for $a\in\hat{\mathcal V}^\DR$. 

One then defines the topological $\mathbf k$-algebra isomorphisms $\mathrm{comp}_{(\mu,g)}^{\mathcal V,(1)}:
\hat{\mathcal V}^\B\to\hat{\mathcal V}^\DR$ and $\mathrm{comp}_{(\mu,g)}^{\mathcal W,(1)}:
\hat{\mathcal W}^\B\to\hat{\mathcal W}^\DR$, and the $\mathbf k$-module isomorphisms 
$\mathrm{comp}_{(\mu,g)}^{\mathcal V,(10)}:\hat{\mathcal V}^\B\to\hat{\mathcal V}^\DR$ and 
$\mathrm{comp}_{(\mu,g)}^{\mathcal M,(10)}:\hat{\mathcal M}^\B\to\hat{\mathcal M}^\DR$ by 
$\mathrm{comp}_{(\mu,g)}^{\mathcal X,(\alpha)}:=\mathrm{aut}_{(\mu,g)}^{\mathcal X,\DR,(\alpha)}
\circ i^{\mathcal X}$ for $\mathcal X\in\{\mathcal V,\mathcal W,\mathcal M\}$ and $\alpha\in\{1,10\}$. 

One then defines ${}^\Gamma\!\mathrm{comp}_{(\mu,g)}^{\mathcal W,(1)}:=\mathrm{Ad}_{\Gamma(\mu,g)}\circ
\mathrm{comp}_{(\mu,g)}^{\mathcal W,(1)}$ and ${}^\Gamma\!\mathrm{comp}_{(\mu,g)}^{\mathcal M,(10)}:=
\ell_{\Gamma(\mu,g)}\circ\mathrm{comp}_{(\mu,g)}^{\mathcal W,(1)}$.

\subsubsection{Geometric interpretations}

In \cite{De}, prounipotent $\mathbb Q$-group schemes $\pi_1^\DR(X,x)$ (\S12.4) and $\pi_1^\B(X,x)$ (\S10.5, denoted there 
$\pi_1^{\text{alg,un}}(X,x)$) are attached to particular schemes $X$ with (tangential) base point $x$; torsors $\pi_1^\DR(X;y,x)$ and $\pi_1^\B(X;y,x)$ are also attached to the datum of an additional (tangential) base point $y$. One has compatible algebra 
and module identifications $\mathrm{ULie}\pi_1^\omega(\mathfrak M_{0,4},\vec 1)\simeq\hat{\mathcal V}_{\mathbb C}^\omega$ 
and $\mathcal O(\pi_1^\omega(\mathfrak M_{0,4},\vec 1,\vec 0))^\vee\simeq\hat{\mathcal V}_{\mathbb C}^\omega$ 
for $\omega\in\{\B,\DR\}$; here $\mathfrak M_{0,4}$ is the moduli space of genus zero curves with four marked points 
and $\vec 0$\index{0 vec@$\vec 0$} and $\vec 1$\index{1 vec@$\vec 1$} are its tangential base points corresponding 
to $(0,{\partial\over\partial z})$ and $(1,-{\partial\over\partial z})$ under the identification 
$\overline{\mathfrak M}_{0,4}\simeq\mathbb P^1$ with coordinate $z$; $\mathrm{ULie}(-)$ is the cocommutative topological 
Hopf algebra attached to a $\mathbb Q$-group scheme and $\mathcal O(-)$ be the ring of regular functions over a scheme. 
When $(\mu,g)=(2\pi i,\varphi_\KZ)$, then $\mathrm{comp}^{\mathcal V,(1)}_{(\mu,g)}$ and 
$\mathrm{comp}^{\mathcal V,(10)}_{(\mu,g)}$ 
can respectively be identified with the Betti-de Rham comparison algebra and module isomorphisms 
$\mathrm{ULie}\pi_1^\B(\mathfrak M_{0,4},\vec 1)\stackrel{\sim}{\to}\mathrm{ULie}\pi_1^\DR(\mathfrak M_{0,4},\vec 1)$ 
and $\mathcal O(\pi_1^\B(\mathfrak M_{0,4},\vec 1,\vec 0))^\vee\stackrel{\sim}{\to}
\mathcal O(\pi_1^\DR(\mathfrak M_{0,4},\vec 1,\vec 0))^\vee$ given by \cite{De}, \S12.16.

\section{Geometric interpretation of the harmonic coproducts}\label{sect:2:27032020}

\subsection{The topology of the moduli spaces $\mathfrak M_{0,4}$ and $\mathfrak M_{0,5}$}
\label{sect:2:3:27032020}

For $n\geq4$, let $\mathfrak M_{0,n}$ be the moduli space of complex curves of genus $0$ with $n$ marked points. 
A contractible subspace of $\mathfrak M_{0,n}$ is the space $b_n$ of marked curves given by $\mathbb P^1(\mathbb C)$ 
with marked points $(x_1,\ldots,x_n)$ distributed on the subset $\mathbb P^1(\mathbb R)$ in counterclockwise order. 
The corresponding fundamental group $P_n^*:=\pi_1(\mathfrak M_{0,n},b_n)$ is the pure modular group of the sphere with $n$ 
marked points; one has $P_n^*\simeq K_{n-1}/Z(K_{n-1})$, where $K_{n-1}$ is the Artin pure braid group of $n-1$ strands on the 
plane, and $Z(K_{n-1})$ is its center, isomorphic to $\mathbb Z$. Recall $K_n$ is the kernel of the natural morphism $B_n\to S_n$, 
where $B_n$ is the Artin braid group with $n$ strands; let $\sigma_1,\ldots,\sigma_{n-1}$ be the Artin generators of $B_n$, 
satisfying the relations $\sigma_i\sigma_{i+1}\sigma_i=\sigma_{i+1}\sigma_i\sigma_{i+1}$ for $i<n-1$ and 
$\sigma_i\sigma_j=\sigma_j\sigma_i$ for $|i-j|\geq2$. Define the family of elements $x_{ij}:=\sigma_{j-1}\cdots\sigma_{i+1}
\sigma_i^2\sigma_{i+1}^{-1}\cdots\sigma_{j-1}^{-1}\in B_n$ for $1\leq i<j\leq n$; one shows that $x_{ij}\in K_n$, and that 
this family generates $K_n$. The center $Z(K_n)$ is generated by 
$x_{12}\cdot x_{13}x_{23}\cdot \cdots\cdot x_{1n}\cdots x_{n-1,n}$. 

For $i\in[1,5]$, let $\underline{\mathrm{pr}}_i:\mathfrak M_{0,5}\to\mathfrak M_{0,4}$ be the map corresponding to the 
erasing of the point labeled $i$, and denote in the same way the induced morphism $P^*_5\to P^*_4$. 
The operation of replacing the point labeled $4$ with two nearby points labeled  $4$ and $5$ induces a morphism 
$\underline\ell:P_4^*\to P_5^*$. We also denote by $\underline{\mathrm{pr}}_i$, $\underline\ell$ the morphisms 
between group algebras induced by these group morphisms. 

When $n=4$, then $P_4^*\simeq K_3/Z(K_3)$ is freely generated by $x_{12}$ and $x_{23}$; we identify it with the group 
$F_2=\mathcal G(\mathcal V^\B)$ of group-like elements of $(\mathcal V^\B,\Delta^{\mathcal V,\B})$ via 
$X_0\simeq x_{23}$, $X_1\simeq x_{12}$. 

To any group $\Gamma$, one functorially associates the $\mathbb Z_{\geq0}$-graded $\mathbb Z$-Lie algebra 
$\mathrm{gr}(\Gamma)$ attached to its lower central series, and therefore the $\mathbb Z_{\geq0}$-graded 
$\mathbf k$-Lie algebra $\mathrm{gr}(\Gamma)\otimes\mathbf k$.  We set $\mathfrak p_n:=\mathrm{gr}(P_n^*)
\otimes\mathbf k$. The Lie algebra $\mathfrak p_n$ is presented by generators $e_{ij}$, $1\leq i\neq j\leq n$ subject to 
relations $e_{ji}=e_{ij}$, $[e_{ij},e_{kl}]=0$ for $i,j,k,l$ distinct, $\sum_{j\in[1,n]-\{i\}}e_{ij}=0$; 
if $i<j$, then $e_{ij}$ is the image of $x_{ij}-1$  in $\mathrm{gr}_1(P_n^*)\otimes\mathbf k\subset\mathfrak p_n$. 
The graded Lie algebra morphisms induced by $\underline\ell$ and $\underline{\mathrm{pr}}_i$ will be denoted
$\ell:\mathfrak p_4\to\mathfrak p_5$ and $\mathrm{pr}_i:\mathfrak p_5\to\mathfrak p_4$. The corresponding 
bialgebra morphisms between universal enveloping algebras will be denoted in the same way. 

When $n=4$, then $\mathfrak p_4$ is freely generated by $e_{12}$ and $e_{23}$; it can be identified with the Lie algebra
$\mathfrak f_2=\mathcal P(\mathcal V^\DR)$ of primitive elements of $(\mathcal V^\DR,\Delta^\DR)$ via $e_0\simeq e_{23}$, 
$e_1\simeq e_{12}$. The isomorphism $\mathrm{gr}(F_2)\otimes\mathbf k\simeq\mathfrak f_2$ takes the image of $X_i-1$ in 
$\mathrm{gr}_1(F_2)\otimes\mathbf k\hookrightarrow \mathfrak f_2$ to $e_i$ for $i=0,1$.  

\subsection{Geometric interpretation of the Betti harmonic coproducts $\Delta^{\mathcal X,\B}$, $\mathcal X
\in\{\mathcal W,\mathcal M\}$}

\subsubsection{Interpretation of $\Delta^{\mathcal W,\B}$}

Let $J(\underline{\mathrm{pr}}_5)$ be the kernel of the algebra morphism $\underline{\mathrm{pr}}_5:\mathbf kP_5^*
\to\mathbf kP_4^*\simeq\mathcal V^\B$. Then  $J(\underline{\mathrm{pr}}_5)$  is a two-sided ideal of $\mathbf kP_5^*$, 
freely generated by the family $(x_{i5}-1)_{1\leq i\leq 3}$ both as a left and as a right $\mathbf kP_5^*$-module. 
These properties imply that for any $a\in\mathbf kP_5^*$, there is a unique element $\underline\varpi(a)=
(\underline\varpi(a)_{ij})_{1\leq i,j\leq 3}\in M_3(\mathbf kP_5^*)$, such that 
$(x_{i5}-1)a=\sum_{i=1}^3\underline\varpi(a)_{ij}(x_{j5}-1)$ for any $i\in[1,3]$, and that the map 
$\underline\varpi:\mathbf kP_5^*\to M_3(\mathbf kP_5^*)$ is an algebra morphism.  
 
Let $\mathcal V^\B[(X_1-1)^{-1}]$ be the localization of $\mathcal V^\B$ with respect to $X_1-1$.
It admits an algebra $\mathbb Z$-filtration given $F^n\mathcal V^\B[(X_1-1)^{-1}]=\sum_{k\geq 0,n_0,\ldots,n_k|
n_0+\cdots+n_k-k=n}F^{n_0}\mathcal V^\B(X_1-1)^{-1}\cdots(X_1-1)^{-1}F^{n_k}\mathcal V^\B$ for $n\in\mathbb Z$.    
Define the elements ${\mathrm{row}}^\B_1:=\begin{pmatrix} 1\otimes(1-X_1^{-1})^{-1}& (1-X_1)^{-1}\otimes1
& 0\end{pmatrix}\in 
M_{1\times 3}(\mathcal V^\B[(X_1-1)^{-1}]^{\otimes2})$, 
${\mathrm{col}}_1^\B:=\begin{pmatrix} (X_1-1)\otimes (X_1-1) \\ (X_1-1)\otimes (X_1^{-1}-1) \\ 0\end{pmatrix}\in 
M_{3\times 1}((\mathcal V^\B)^{\otimes2})$.  

\begin{prop} \label{prop:2:1} (see \cite{EF1}, Proposition 8.6) The diagram 
$$
\xymatrix{
\mathcal V^\B\ar^{\underline\ell}[r]\ar^\sim_{(-)\cdot(X_1-1)}[d]& \mathbf kP_5^*\ar^{\!\!\!\underline\varpi}[r]& 
M_3(\mathbf kP_5^*)\ar^{\!\!M_3(\underline{\mathrm{pr}}_{12})}[r]& M_3((\mathcal V^\B)^{\otimes2})
\ar^{\!\!\!\!{\mathrm{row}}_1^\B\cdot(-)\cdot {\mathrm{col}}_1^\B}[rr]&&
\mathcal V^\B[(X_1-1)^{-1}]^{\otimes2} \\
\mathcal W^\B_+\ar_{\Delta^{\mathcal W,\B}}[rrrrr]&& & & & (\mathcal W^\B)^{\otimes2}\ar@{^{(}->}[u]}
$$
is commutative, where $\underline{\mathrm{pr}}_{12}:\mathbf kP_5^*\to(\mathcal V^\B)^{\otimes2}$ is the morphism induced
by the group morphism $P_5^*
\stackrel{\mathrm{diag}}{\longrightarrow}
(P_5^*)^2
\stackrel{\underline{\mathrm{pr}}_1\times\underline{\mathrm{pr}}_2}{\longrightarrow}
(P_4^*)^2$, where $\mathrm{diag}$ is the diagonal morphism, where 
${\mathrm{row}}_1^\B\cdot(-)\cdot {\mathrm{col}}^\B_1$ is the map 
$a\mapsto{\mathrm{row}}^\B_1\cdot a\cdot {\mathrm{col}}^\B_1$, where $\mathcal W^\B_+$ is the subalgebra without 
unit $\mathcal V^\B(X_1-1)\hookrightarrow\mathcal W^\B$, and the left vertical map is $a\mapsto a(X_1-1)$. In this 
diagram, all the maps are compatible with the filtrations, except for the left vertical and the rightmost top horizontal 
maps, which increase the filtration degree by one. 
\end{prop}

{\it Idea of proof.} The proof is based on: (1) the fact that $\underline\rho:=M_3(\underline{\mathrm{pr}}_{12})\circ
\underline\varpi\circ\underline\ell$ is an algebra morphism such that $X_1-1\mapsto{\mathrm{col}}_1^\B\cdot
{\mathrm{row}}_1^\B$; (2) the fact that (1) implies that $({\mathrm{row}}^\B_1\cdot(-)\cdot
{\mathrm{col}}^\B_1)\circ\underline\rho\circ((-)\cdot(X_1-1))$ is an algebra morphism $\mathcal W_\B^+\to
\mathcal V^\B[(X_1-1)^{-1}]^{\otimes2}$; (3) the identification of a generating family of $\mathcal W_\B^+$; (4) the 
identification of the images by $({\mathrm{row}}^\B_1\cdot(-)\cdot{\mathrm{col}}^\B_1)\circ
\underline\rho\circ((-)\cdot(X_1-1))$ and $\Delta^{\mathcal W,\B}$ of each element of the generating family of (3). 
\hfill\qed\medskip 

\subsubsection{Interpretation of $\Delta^{\mathcal M,\B}$}

Set $\mathcal M^\B[(X_1-1)^{-1}]:=\mathcal V^\B[(X_1-1)^{-1}]/\mathcal V^\B[(X_1-1)^{-1}](X_0-1)$. 
This is a left $\mathcal V^\B[(X_1-1)^{-1}]$-module. There is a natural morphism $\mathcal M^\B\to
\mathcal M^\B[(X_1-1)^{-1}]$ compatible with the algebra morphism $\mathcal V^\B\to\mathcal V^\B[(X_1-1)^{-1}]$, 
and which can be shown to be injective. We denote by $1_\B\in\mathcal M^\B[(X_1-1)^{-1}]$ the image of 
$1_\B\in\mathcal M^\B$. The module $\mathcal M^\B[(X_1-1)^{-1}]$ is equipped with a $\mathbb Z$-filtration 
given by $F^n\mathcal M^\B[(X_1-1)^{-1}]:=F^n\mathcal V^\B[(X_1-1)^{-1}]\cdot 1_\B$ for $n\in\mathbb Z$
and is then a filtered module over $\mathcal V^\B[(X_1-1)^{-1}]$. 
Let ${\mathrm{col}}_0^\B:=\begin{pmatrix} 0\\ ((1-X_1)\otimes X_1^{-1})\cdot 1_\B^{\otimes2} \\ 
((1-X_1^{-1})\otimes X_1^{-1})\cdot 1_\B^{\otimes2}\end{pmatrix}\in 
M_{3\times 1}((\mathcal M^\B)^{\otimes2})$.

\begin{prop} \label{prop:2:2} (see \cite{EF1}, Proposition 8.14) The diagram 
$$
\xymatrix{
\mathcal V^\B\ar^{\underline\ell}[r]\ar_{(-)\cdot 1_\B}[d]& \mathbf kP_5^*\ar^{\!\!\!\underline\varpi}[r]& 
M_3(\mathbf kP_5^*)\ar^{\!\!M_3(\underline{\mathrm{pr}}_{12})}[r]& M_3((\mathcal V^\B)^{\otimes2})
\ar^{\!\!\!\!{\mathrm{row}}_1^\B\cdot(-)\cdot {\mathrm{col}}_0^\B}[rr]&&
\mathcal M^\B[(X_1-1)^{-1}]^{\otimes2} \\
\mathcal M^\B\ar_{\Delta^{\mathcal M,\B}}[rrrrr]&& & & & (\mathcal M^\B)^{\otimes2}\ar@{^{(}->}[u]}
$$
is commutative, where ${\mathrm{row}}_1^\B\cdot(-)\cdot {\mathrm{col}}^\B_0$ is the map 
$a\mapsto{\mathrm{row}}^\B_1\cdot a\cdot {\mathrm{col}}^\B_0$. In this 
diagram, all the maps are compatible with the filtrations. 
\end{prop}

{\it Idea of proof.} The proof is based on: (1) the fact that $\underline\rho(X_0-1)\cdot\mathrm{col}_0^\B=0$, 
which implies that there exists a map $\delta:\mathcal M^\B\to\mathcal M^\B[(X_1-1)^{-1}]^{\otimes2}$, such that 
$\delta\circ((-)\cdot 1_\B)=(\mathrm{row}_1^\B\cdot(-)\cdot\mathrm{col}_0^\B)\circ\underline\rho$ ($\underline\rho$ 
being as in the proof of Proposition \eqref{prop:2:1}); (2) Proposition \eqref{prop:2:1}, which implies that $\delta$ is 
compatible with the module structures on both sides and with the morphism $\Delta^{\mathcal W,\B}:\mathcal W^\B\to
(\mathcal W^\B)^{\otimes2}$; (3) the fact that $\mathrm{row}_1^\B\cdot\mathrm{col}_0^\B=1_\B^{\otimes2}$, which 
implies  that $\delta(1_\B)=1_\B^{\otimes2}$ by (1), and therefore $\delta=\Delta^{\mathcal W,\B}$ by (2). 
\hfill\qed\medskip 

\subsection{Geometric interpretation of the de Rham harmonic coproducts $\Delta^{\mathcal X,\DR}$, 
$\mathcal X\in\{\mathcal W,\mathcal M\}$}

\subsubsection{Interpretation of $\Delta^{\mathcal W,\DR}$}

The morphism $\underline\varpi:\mathbf kP_5^*\to M_3(\mathbf kP_5^*)$ is compatible with the filtrations. Let us denote by 
$\varpi:U(\mathfrak p_5)\to M_3(U(\mathfrak p_5))$ the associated graded morphism. One can show that $\varpi$ may be constructed as follows. Let $J(\mathrm{pr}_5)$ be the kernel of the algebra morphism $\mathrm{pr}_5:U(\mathfrak p_5)
\to U(\mathfrak p_4)$; this is two-sided ideal of $U(\mathfrak p_5)$, free and generated by $(e_{i5})_{1\leq i\leq 3}$
both as a left and right $U(\mathfrak p_5)$-module. For $a\in U(\mathfrak p_5)$, the matrix 
$\varpi(a)=(\varpi(a)_{ij})_{1\leq i,j\leq 3}$ is uniquely determined by the identity $e_{i5}a=\sum_{i=1}^3 
\varpi(a)_{ij}e_{j5}$ for $1\leq i\leq 3$. 

Let $\mathcal V^\DR[e_1^{-1}]$ be the localization of $\mathcal V^\DR$ with respect to $e_1$. This is a 
$\mathbb Z$-graded algebra, which can be identified with $\mathrm{gr}(\mathcal V^\B[(X_1-1)^{-1}])$. 
Define the elements $\mathrm{row}_1^\DR:=\begin{pmatrix} 1\otimes e_1^{-1}& -e_1^{-1}\otimes1& 0\end{pmatrix}\in 
M_{1\times3}(\mathcal V^\DR[e_1^{-1}]^{\otimes2})$ and $\mathrm{col}_1^\DR:=\begin{pmatrix} e_1\otimes e_1\\ 
-e_1\otimes e_1\\ 0  \end{pmatrix}\in 
M_{1\times3}(\mathcal V^\DR[e_1^{-1}]^{\otimes2})$. 

\begin{prop} \label{prop:2:3} (see \cite{DeT} and \cite{EF1}, Proposition 6.3) The diagram 
$$
\xymatrix{
\mathcal V^\DR\ar^{\ell}[r]\ar^\sim_{(-)\cdot e_1}[d]& U(\mathfrak p_5)\ar^{\!\!\!\varpi}[r]& 
M_3(U(\mathfrak p_5))\ar^{\!\!M_3({\mathrm{pr}}_{12})}[r]& M_3((\mathcal V^\DR)^{\otimes2})
\ar^{\!\!\!\!{\mathrm{row}}_1^\DR\cdot(-)\cdot {\mathrm{col}}_1^\DR}[rr]&&
\mathcal V^\DR[e_1^{-1}]^{\otimes2} \\
\mathcal W^\DR_+\ar_{\Delta^{\mathcal W,\DR}}[rrrrr]&& & & & (\mathcal W^\DR)^{\otimes2}\ar@{^{(}->}[u]}
$$
is commutative, where ${\mathrm{pr}}_{12}:U(\mathfrak p_5)\to(\mathcal V^\DR)^{\otimes2}$ is the morphism induced
by the group morphism $\mathfrak p_5\stackrel{\mathrm{diag}}{\longrightarrow}(\mathfrak p_5)^2
\stackrel{{\mathrm{pr}}_1\times{\mathrm{pr}}_2}{\longrightarrow}
(\mathfrak p_4)^2$, where $\mathrm{diag}$ is the diagonal morphism, where 
${\mathrm{row}}_1^\DR\cdot(-)\cdot {\mathrm{col}}^\DR_1$ is the map 
$a\mapsto{\mathrm{row}}^\DR_1\cdot a\cdot {\mathrm{col}}^\DR_1$, where $\mathcal W^\DR_+$ is the subalgebra without 
unit $\mathcal V^\DR e_1\hookrightarrow\mathcal W^\DR$, and the left vertical map is $a\mapsto ae_1$. In this 
diagram, all the maps are of degree zero, except for the left vertical and the rightmost top horizontal 
maps, which are of degree one. 
\end{prop}

{\it Idea of proof.} This can be proved either by repeating the steps of the proof of Proposition \ref{prop:2:1}, 
which is what is being done in \cite{EF1}, or by applying the associated graded functor to the diagram of this 
proposition. \hfill\qed\medskip 

\subsubsection{Interpretation of $\Delta^{\mathcal M,\DR}$}

Set $\mathcal M^\DR[e_1^{-1}]:=\mathcal V^\DR[e_1^{-1}]/\mathcal V^\DR[e_1^{-1}]e_0$. 
This is a left $\mathcal V^\DR[e_1^{-1}]$-module. There is a natural morphism $\mathcal M^\DR\to
\mathcal M^\DR[e_1^{-1}]$ compatible with the algebra morphism $\mathcal V^\DR\to\mathcal V^\DR[e_1^{-1}]$, 
and which can be shown to be injective. We denote by $1_\DR\in\mathcal M^\DR[e_1^{-1}]$ the image of 
$1_\DR\in\mathcal M^\DR$. The module $\mathcal M^\DR[e_1^{-1}]$ is equipped with a $\mathbb Z$-grading
defined by the condition that $1_\DR$ has degree 0 and that $\mathcal M^\DR[e_1^{-1}]$ is graded as a 
$\mathcal V^\DR[e_1^{-1}]$-module. 
Let ${\mathrm{col}}_0^\DR:=\begin{pmatrix} 0\\ -(e_1\otimes1)\cdot 1_\DR^{\otimes2} \\ 
(e_1\otimes 1)\cdot 1_\DR^{\otimes2}\end{pmatrix}\in 
M_{3\times 1}((\mathcal M^\DR)^{\otimes2})$.

\begin{prop} \label{prop:2:4} (see \cite{EF1}, Proposition 6.9) The diagram 
$$
\xymatrix{
\mathcal V^\DR\ar^{\ell}[r]\ar_{(-)\cdot 1_\DR}[d]& U(\mathfrak p_5)\ar^{\!\!\!\varpi}[r]& 
M_3(U(\mathfrak p_5))\ar^{\!\!M_3({\mathrm{pr}}_{12})}[r]& M_3((\mathcal V^\DR)^{\otimes2})
\ar^{\!\!\!\!{\mathrm{row}}_1^\DR\cdot(-)\cdot {\mathrm{col}}_0^\DR}[rr]&&
\mathcal M^\DR[e_1^{-1}]^{\otimes2} \\
\mathcal M^\DR\ar_{\Delta^{\mathcal M,\DR}}[rrrrr]&& & & & (\mathcal M^\DR)^{\otimes2}\ar@{^{(}->}[u]}
$$
is commutative, where ${\mathrm{row}}_1^\DR\cdot(-)\cdot {\mathrm{col}}^\DR_0$ is  
$a\mapsto{\mathrm{row}}^\DR_1\cdot a\cdot {\mathrm{col}}^\DR_0$. In this 
diagram, all the maps have degree zero. 
\end{prop}

{\it Idea of proof.} This can be proved either by repeating the steps of the proof of Proposition \ref{prop:2:2}, 
which is what is being done in \cite{EF1}, or by applying the associated graded functor to the diagram of that 
proposition. \hfill\qed\medskip 

\section{Associators and double shuffle relations}\label{sect:3:27032020}

\subsection{Associators}

The notion of associator was defined in \cite{Dr}; it was then shown in \cite{F2} that this definition can 
be formulated as follows. For $(I,J,K)$ one of the triples $(2,3,4)$, $(12,3,4)$, $(1,23,4)$, $(1,2,34)$, $(1,2,3)$, 
define a Lie algebra morphism $\mathfrak f_2\to\mathfrak p_5$, $x\mapsto x^{I,J,K}$ by $e_0^{I,J,K}:=e_{I,J}$ 
and $e_1^{I,J,K}:=e_{J,K}$, where $e_{ij,k}:=e_{ik}+e_{jk}$ and $e_{i,jk}:=e_{jk,i}$. We denote the induced 
algebra morphisms $\hat{\mathcal V}^\DR\to U(\mathfrak p_5)^\wedge$ in the same way. 

\begin{defn}
For $\mu\in\mathbf k^\times$, one sets $M_\mu(\mathbf k):=\{\Phi\in\mathcal G(\hat{\mathcal V}^\DR)|
\Phi^{2,3,4}\Phi^{1,23,4}\Phi^{1,2,3}=\Phi^{12,3,4}\Phi^{1,2,34}$, $(\Phi|e_0)=(\Phi|e_1)=0$ and 
$(\Phi|e_0e_1)={\mu^2\over{24}}\}$, where $(\Phi|e_0e_1)$
is the coefficient of $e_0e_1$ as an expansion in words in $e_0,e_1$. 
The set of associators is $M(\mathbf k):=\sqcup_{\mu\in\mathbf k^\times}M_\mu(\mathbf k)$. 
\end{defn}

One has $\varphi_{\mathrm{KZ}}\in M_{2\pi i}(\mathbb C)$.

\subsection{Compatibility of the associators with the coproducts}

\begin{thm} \label{thm:part:A:2804}(\cite{EF1}, Theorems 10.9 and 11.13)
Let $\mu\in\mathbf k^\times$, and $\Phi\in M_\mu(\mathbf k)$. 
 
1) The diagram of $\mathbf k$-algebra morphisms
$$
\xymatrix{
\hat{\mathcal W}^\B\ar^{\hat\Delta^{\mathcal W,\B}}[r]\ar_{^\Gamma\!\mathrm{comp}_{(\mu,\Phi)}^{\mathcal W,(1)}}[d] 
& (\mathcal W^\B)^{\otimes2\wedge}\ar^{(^\Gamma\!\mathrm{comp}_{(\mu,\Phi)}^{\mathcal W,(1)})^{\otimes2}}[d] \\ 
\hat{\mathcal W}^\DR\ar_{\hat\Delta^{\mathcal W,\DR}}[r] 
& (\mathcal W^\DR)^{\otimes2\wedge}}
$$
is commutative. 

2) The diagram of $\mathbf k$-module morphisms
$$
\xymatrix{
\hat{\mathcal M}^\B\ar^{\hat\Delta^{\mathcal M,\B}}[r]\ar_{^\Gamma\!\mathrm{comp}_{(\mu,\Phi)}^{\mathcal M,(10)}}[d] 
& (\mathcal M^\B)^{\otimes2\wedge}\ar^{(^\Gamma\!\mathrm{comp}_{(\mu,\Phi)}^{\mathcal M,(10)})^{\otimes2}}[d] \\ 
\hat{\mathcal M}^\DR\ar_{\hat\Delta^{\mathcal M,\DR}}[r] 
& (\mathcal M^\DR)^{\otimes2\wedge}}
$$
is commutative. 
\end{thm}

{\it Sketch of proof of 1).} In \cite{BN}, one introduces the categories $\mathbf{PaB}$ 
and $\mathbf{PaCD}$ of parenthesized braids and parenthesized braid 
diagrams, whose sets of objects both coincide with the set of parenthesized 
words in one letter $\bullet$, and one attaches to each associator 
$(\mu,\Phi)$ a functor $\mathbf{PaB}\to\mathbf{PaCD}$. If $P,Q$ are two 
parenthesized words in $\bullet$ of length $n$, and $\mathcal C(X,Y)$ is the set of morphisms 
$X\to Y$, where $X,Y$ are objects of a category $\mathcal C$, ones denotes by  
$\mathrm{comp}_{(\mu,\Phi)}^{P,Q}:B_n\simeq\mathbf{PaB}(P,Q)
\to\mathbf{PaCD}(P,Q)\simeq U(\mathfrak p_{n+1})^\wedge\rtimes S_n$ 
the resulting map; we also denote $\mathrm{comp}_{(\mu,\Phi)}^{P}$ 
(resp. $\mathbf{PaB}(P)$, $\mathbf{PaCD}(P)$) for 
$\mathrm{comp}_{(\mu,\Phi)}^{P,P}$ (resp. 
$\mathbf{PaB}(P,P)$, $\mathbf{PaCD}(P,P)$). 

Define the element 
$P_{(\mu,\Phi)}\in M_3((U\mathfrak p_5)^{\otimes 2\wedge})$ by the equality 
$$
M_{3\times 1}(\mathrm{comp}_{(\mu,\Phi)}^{((\bullet\bullet)\bullet)\bullet})
\begin{pmatrix} x_{15}-1 \\ x_{25}-1 \\ x_{35}-1\end{pmatrix}
=P_{(\mu,\Phi)}\begin{pmatrix} e_{15} \\ e_{25} \\ e_{35}\end{pmatrix}. 
$$
and set 
$$
\overline P_{(\mu,\Phi)}:=M_3(\mathrm{pr}_{12})(P_{(\mu,\Phi)})
\in M_3((\mathcal V^\DR)^{\otimes 2\wedge}). 
$$
In \cite{EF1}, Proposition 9.20 and Corollary 9.21, one gives 
explicit expressions for $P_{(\mu,\Phi)}$ and $\overline P_{(\mu,\Phi)}$. 

Define also the following elements in $((\mathcal W^\DR)^{\otimes2\wedge})^\times\subset ((\mathcal V^\DR)^{\otimes2\wedge})^\times$:  
$$
B_\Phi:={{\Gamma_\Phi(-e_1)\otimes\Gamma_\Phi(-e_1)}\over{
\Gamma_\Phi(-e_1\otimes1-1\otimes e_1)}}, \quad 
u_{(\mu,\Phi)}:=\mu(1\otimes e^{-\mu e_1})
{{e^{\mu(e_1\otimes1+1\otimes e_1)}-1}\over{e_1\otimes1+1\otimes e_1}}
{{\Gamma_\Phi(e_1\otimes1+1\otimes e_1)}\over
{\Gamma_\Phi(e_1)\otimes\Gamma_\Phi(e_1)}}
$$
$$
v_{(\mu,\Phi)}:={1\over\mu}(1\otimes e^{\mu e_1})
{{\Gamma_\Phi(e_1)\otimes\Gamma_\Phi(e_1)}
\over{\Gamma_\Phi(e_1\otimes1+1\otimes e_1)}}, \quad 
\kappa_{(\mu,\Phi)}:=\Phi(e_0,e_1)\otimes
(e^{-(\mu/2)e_1}\Phi(e_\infty,e_1))
$$

Set also 
$$
\underline{\mathrm{row}}_1:=((X_1-1)\otimes(1-X_1^{-1}))
\cdot \mathrm{row}_1^\B\in M_{1\times3}((\mathcal V^\B)^{\otimes2\wedge}), 
$$
$$
\underline{\mathrm{col}}_1:=\mathrm{col}_1^\B\cdot((X_1-1)\otimes(1-X_1^{-1}))^{-1}
\in M_{3\times1}((\mathcal V^\B)^{\otimes2\wedge}); 
$$
$$\mathrm{row}_1:=(e_1\otimes e_1)\cdot
\mathrm{row}_1^\DR\in M_{1\times3}((\mathcal V^\DR)^{\otimes2\wedge}), 
\quad \mathrm{col}_1:=\mathrm{col}_1^\DR\cdot (e_1\otimes e_1)^{-1}
\in M_{3\times1}((\mathcal V^\DR)^{\otimes2\wedge}). 
$$

Define the following morphisms of topological $\mathbf k$-algebras: 
\begin{itemize}
\item $\hat\rho:=(M_3(\mathrm{pr}_{12})\circ\varpi\circ\ell)^\wedge :
\hat{\mathcal V}^\DR\to M_3((\mathcal V^\DR)^{\otimes2})$,
$\underline{\hat\rho}:=(M_3(\underline{\mathrm{pr}}_{12})\circ
\underline\varpi\circ\underline\ell)^\wedge : \hat{\mathcal V}^\B
\to M_3((\mathcal V^\B)^{\otimes2})$, 
where $(-)^\wedge$ indicates the completion with respect to the 
underlying filtrations, 

\item $\mathrm{compl}^{\mathcal V_{\mathrm{loc}},(1)}_{(\mu,\Phi)}$
is the extension of $\mathrm{compl}^{\mathcal V,(1)}_{(\mu,\Phi)}$
to a $\mathbf k$-algebra isomorphism  $\hat{\mathcal V}^\B[(X_1-1)^{-1}]\to
\hat{\mathcal V}^\DR[e_1^{-1}]$,  

\item for $A$ a $\mathbf k$-algebra and $a\in A$, 
$\cdot_a$ is the product on $A$ defined by $x\cdot_a y=xay$, and 
$\mathrm{mor}_{A,a}:(A,\cdot_a)\to A$ is the algebra morphism given by 
$x\mapsto xa$. 
\end{itemize}
Consider the diagram 
$$
{\tiny\xymatrix{
\hat{\mathcal W}^\B_+
\ar^{{\hat\Delta}^{\mathcal W,\B}}[rrrrrr]\ar^*[@]{\hbox to 15pt {$\mathrm{comp}_{(\mu,\Phi)}^{\mathcal W,(1)}$}}[ddddd]
\ar@{}^*[@]{{\hbox to -10pt{$\!\!\!\!\!\!\!\!\!\!\!\!\!\!\mathrm{mor}_{\hat{\mathcal V}^\B,X_1-1}$}}}[rd] 
\ar@{}[dddddr]|{\mathrm{(A2)}}
&&&&& 
\ar@{}[dllll]|{\mathrm{(A1)}}
& 
(\mathcal W^\B)^{\otimes2\wedge}
\ar^*[@]{\hbox to 15pt {$(\mathrm{comp}_{(\mu,\Phi)}^{\mathcal W,(1)})^{\otimes2}$}}[ddd]\ar@{^(->}[ld] \\
& 
\hat{\mathcal V}^\B
\ar^{\hbox to 20pt{$\underline{\hat\rho}$}}[r]
\ar^*[@]{\vbox to 4pt{\hbox to 4pt {$\!\!\!\!\!\!\mathrm{comp}_{(\mu,\Phi)}^{\mathcal V,(1)}$}}}[dd]
\ar[lu]& 
{\scriptstyle M_3((\mathcal V^\B)^{\otimes2\wedge})}
\ar^{\ \ \ \ \underline{\mathrm{row}}_1\cdot(-)\cdot\underline{\mathrm{col}}_1}[rr]
\ar^{M_3((\mathrm{comp}_{(\mu,\Phi)}^{\mathcal V,(1)})^{\otimes2})}[d] && 
(\mathcal V^\B)^{\otimes2\wedge}
\ar^{\begin{matrix} \mathrm{Ad}((X_1-1)^{-1}\cdot\\ \cdot(1-Y_1^{-1})^{-1})\end{matrix}}[r]
\ar^{(\mathrm{comp}_{(\mu,\Phi)}^{\mathcal V,(1)})^{\otimes2}}[d]
&  
F^0\mathcal V^\B[{1\over{X_1-1}}]^{\otimes2\wedge}
\ar^{(\mathrm{comp}_{(\mu,\Phi)}^{\mathcal V_{\mathrm{loc}},(1)})^{\otimes2}}[d] & \\
&\ar@{}[r]|{\mathrm{\!\!\!\!\!\!(A3)}}& 
{\scriptstyle M_3((\mathcal V^\DR)^{\otimes2\wedge})} 
\ar^{\mathrm{Ad}((\kappa_{(\mu,\Phi)}\cdot\overline P_{(\mu,\Phi)})^{-1})}[d]
\ar@{}[rr]|{\mathrm{(A4)}}
&& 
(\mathcal V^\DR)^{\otimes2\wedge}
\ar^{u_{(\mu,\Phi)}^{-1}\cdot(-)\cdot v_{(\mu,\Phi)}^{-1}}[d]
\ar@{}[r]|{\mathrm{(A5)}}
& 
(\mathcal V^\DR[{1\over{e_1}}]_{\geq0})^{\otimes2\wedge} 
\ar^{\begin{matrix} B_\Phi^{-1}\cdot(-)\cdot B_\Phi\cdot \\ \cdot {{e_1+f_1}\over{{e^{\mu(e_1+f_1)}}-1}}\end{matrix}}[d] 
& \\
& 
\hat{\mathcal V}^\DR
\ar_{\hbox to 20pt{$\hat\rho$}}[r]\ar^{\mathrm{mor}_{\hat{\mathcal V}^\DR,e_1}}[d] & 
{\scriptstyle M_3((\mathcal V^\DR)^{\otimes2\wedge})}
\ar_{\mathrm{row}_1\cdot(-)\cdot\mathrm{col}_1}[rr] && 
(\mathcal V^\DR)^{\otimes2\wedge}
\ar_{\mathrm{Ad}_{(e_1f_1)^{-1}}}[r] 
\ar@{}[dll]|{\mathrm{(A7)}} & 
(\mathcal V^\DR[{1\over e_1}]_{\geq0})^{\otimes2\wedge} 
\ar@{}[dr]|{\mathrm{(A6)}}
&
({\mathcal W}^{\DR})^{\otimes2\wedge}
\ar^*[@]{\hbox to 10pt{$\mathrm{Ad}(B_\Phi^{-1})$}}[dd] 
\\
& 
\hat{\mathcal W}^{\DR}_+
\ar_{\hat\Delta^{\mathcal W,\DR}}[rrrr]
\ar[ld] 
&&&
\ar@{}[dll]|{\mathrm{(A8)}}
& 
(\mathcal W^\DR)^{\otimes2\wedge}
\ar@{^(->}[u]
\ar@{}^*[@]{\hbox to 23pt{\vbox to 0pt {$\!\!\!\!\!\!\!\!\!\!\!\!\!\!\!\!\!\!\!\!\!\!\!(-)\cdot{{e_1+f_1}\over{e^{\mu(e_1+f_1)}-1}}$}}}[rd]  
& \\
\hat{\mathcal W}^{\DR}_+
\ar_{\hat\Delta^{\mathcal W,\DR}}[rrrrrr]
\ar@{}^*[@]{{\hbox to 5pt{$\!\!\!\!\!\!\!\!\!\!\!(-)\cdot{e^{\mu e_1}-1}\over{e_1}$}}}[ru] 
&&&&&&(\mathcal W^{\DR})^{\otimes 2\wedge}\ar[ul]
}}
$$
where we write $\mathrm{Ad}(a)$ for $\mathrm{Ad}_a$. In this diagram, the commutativities of (A1) and (A7)
follows from Proposition \ref{prop:2:1} and \ref{prop:2:3}. The commutativities of (A2), (A5) , (A6) and (A8) are 
immediate. The commutativity of (A3) follows from the definition of $\overline P_{(\mu,\Phi)}$, and more precisely
from the relations between $P_{(\mu,\Phi)}$, $\mathrm{comp}_{(\mu,\Phi)}^{((\bullet\bullet)\bullet)\bullet}$ and 
$\hat\varpi$, $\underline{\hat\varpi}$. The commutativity of (A4) is a consequence of two equalities, one in 
$M_{1\times 3}((\mathcal V^\DR)^{\otimes 2\wedge})$ and the other in 
$M_{3\times 1}((\mathcal V^\DR)^{\otimes 2\wedge})$. Both follow from explicit computation based on 
the already mentioned computation of $\overline P_{(\mu,\Phi)}$. One easily derives the commutativity of the 
announced diagram.  

{\it Sketch of proof of 2).} Set 
$$
R_{(\mu,\Phi)}:=\hat\rho(\Phi)^{-1}\overline P_{(\mu,\Phi)}^{-1}
\kappa_{(\mu,\Phi)}^{-1}(\Phi(e_0,e_1)\otimes\Phi(e_0,e_1))\in
\mathrm{GL}_3((\mathcal V^\DR)^{\otimes2\wedge}).
$$
Define the element $Q_{(\mu,\Phi)}\in 
M_3((U\mathfrak p_5)^{\otimes 2\wedge})$ by the equality 
$$
M_{3\times 1}(\mathrm{comp}_{(\mu,\Phi)}^{(\bullet(\bullet\bullet))\bullet})
\begin{pmatrix} x_{15}-1 \\ x_{25}-1 \\ x_{35}-1\end{pmatrix}
=Q_{(\mu,\Phi)}\begin{pmatrix} e_{15} \\ e_{25} \\ e_{35}\end{pmatrix}. 
$$
and set 
$$
\overline Q_{(\mu,\Phi)}:=M_3(\mathrm{pr}_{12})(Q_{(\mu,\Phi)})
\in M_3((\mathcal V^\DR)^{\otimes 2\wedge}). 
$$
In \cite{EF1}, Proposition 9.23 and Corollary 9.24, one gives 
explicit expressions for $Q_{(\mu,\Phi)}$ and $\overline Q_{(\mu,\Phi)}$. 

By the categorical origin of the morphisms $\mathrm{comp}^P_{(\mu,P)}$,
the two algebra morphisms corresponding to $P=((\bullet\bullet)\bullet)
\bullet$ and $(\bullet(\bullet\bullet))\bullet$ are related by an inner 
conjugation. One derives from there an expression of $Q_{(\mu,\Phi)}$ 
in terms of $P_{(\mu,\Phi)}$ (see \cite{EF1}, Lemma 11.6), which implies 
the relation 
$$
\overline Q_{(\mu,\Phi)}=\overline P_{(\mu,\Phi)}\hat\rho(\Phi)
$$ 
(\cite{EF1}, Lemma 11.7), from which one derives  
\begin{equation}\label{eq:R:mu:Phi}
R_{(\mu,\Phi)}^{-1}=(\Phi(e_0,e_1)\otimes\Phi(e_0,e_1))^{-1}
\kappa_{(\mu,\Phi)}\overline Q_{(\mu,\Phi)}
\end{equation}
(\cite{EF1}, Corollary 11.8). 

Consider the diagram 
$$
{\xymatrix{
{\hat{\mathcal M}^\B} 
\ar^{{\hat\Delta}^{\mathcal M,\B}}[rr]
\ar^*[@]{\vbox to 4pt{\hbox to 4pt {$\!\!\!\!\!\!\scriptstyle{\mathrm{comp}^{\mathcal M,(10)}_{(\mu,\Phi)}}$}}}[dddd]
& & { ({\mathcal M}^\B)^{\otimes2\wedge}}\ar@{^(->}[rr]& & 
F^{-1}({{\mathcal M}^\B[\frac{1}{X_1-1}])^{\otimes 2\wedge}}
\ar^{(\mathrm{comp}^{{\mathcal M}_{\mathrm{loc}},(10)}_{(\mu,\Phi)})^{\otimes2}}[dddd] 
\\
&\hat{\mathcal V}^\B
\ar@{}[ur]|{\mathrm{(M1)}}
\ar@{}[dr]|{\mathrm{(M3)}}
\ar@{}[dddl]|{\mathrm{(M2)}}
\ar^{\!\!\!\!\!\!\!\underline{\hat\BB}}[r]
\ar_{(-)\cdot 1_\B}[ul]
\ar^*[@]{\vbox to 4pt{\hbox to 4pt {$\!\!\!\!\!\!\scriptstyle{\mathrm{comp}^{\mathcal V,(10)}_{(\mu,\Phi)}}$}}}[dd]
&
{ M_3((\mathcal V^\B)^{\otimes 2\wedge})}
\ar^*[@]{\vbox to 4pt{\hbox to 10pt {$\!\!\!\!\!\!\!\!\!\!\!\!\!\!\!\!\!\!\!\!\!\!\!\!\!\!\!\!\!\!\!\!\!\!\!\!\!\scriptstyle{
{\mathrm{row}}_1^\B\cdot(\text{-})\cdot{\mathrm{col}}_0^\B}$}}}[urr]
\ar^{M_3((\mathrm{comp}^{\mathcal V,(10)}_{(\mu,\Phi)})^{\otimes2})}[d]
\ar@{}[ddrr]|{\mathrm{(M4)}}&&
\\
&&
{M_3((\mathcal V^\DR)^{\otimes2\wedge})}
\ar^{(\kappa_{(\mu,\Phi)}\overline P_{(\mu,\Phi)})^{-1}\cdot (\text{-})\cdot R_{(\mu,\Phi)}^{-1}}[d]
&&
\\
& \hat{\mathcal V}^\DR
\ar^{\!\!\!\!\!\!\hat\BB}[r]
\ar[dl]
\ar@{}[dr]|{\mathrm{(M5)}}
& M_3((\mathcal V^\DR)^{\otimes2\wedge})
\ar^*[@]{\vbox to 4pt{\hbox to 10pt {$\!\!\!\!\!\!\!\!\!\!\!\!\!\!\scriptstyle{
\mathrm{row}_1^\DR\cdot(\text{-})\cdot\mathrm{col}_0^\DR}$}}}[rd]
&&
\\
\hat{\mathcal M}^\DR
\ar_{\hat\Delta^{\mathcal M,\DR}}[rr]
\ar@{}^*[@]{{\hbox to 0pt{\vbox to 10pt{$\!\!\!\!\!\!\!\!\!\!\!\!\scriptstyle{(-)\cdot 1_\DR}$}}}}[ur] 
&&({\mathcal M}^\DR)^{\otimes2\wedge}\ar@{^(->}[r]&
({\mathcal M}^\DR[\frac{1}{e_1}]_{\geq-1})^{\otimes2\wedge}
\ar_{\scriptstyle{B_\Phi\cdot(\text{-})}}[r]
&
({\mathcal M}^\DR[\frac{1}{e_1}]_{\geq-1})^{\otimes2\wedge}
}}
$$
where $\mathrm{comp}^{\mathcal M_{\mathrm{loc}},(10)}_{(\mu,\Phi)} : 
\mathcal M^\B[(X_1-1)^{-1}]^\wedge\to\mathcal M^\DR[e_1^{-1}]^\wedge$ 
is the $\mathbf k$-module isomorphism which both extends the 
$\mathbf k$-module isomorphism $\mathrm{comp}^{\mathcal M,(10)}_{(\mu,\Phi)}:
\hat{\mathcal M}^\B\to\hat{\mathcal M}^\DR$ and is compatible with the 
$\mathbf k$-algebra isomorphism
$\mathrm{comp}^{\mathcal V_{\mathrm{loc}},(1)}_{(\mu,\Phi)} : 
\mathcal V^\B[(X_1-1)^{-1}]^\wedge\to\mathcal V^\DR[e_1^{-1}]^\wedge$.  

In this diagram, the commutativity of (M2) is immediate. The commutativities of (M1) and (M5)
are consequences of Proposition \ref{prop:2:2} and \ref{prop:2:4}. 

(M3) states the equality of two $\mathbf k$-module morphisms 
$\hat{\mathcal M}^\B\to M_3((\mathcal M^\DR)^{\otimes2\wedge})$, which 
turn out to be free rank one module morphisms over the two algebra 
morphisms $\hat{\mathcal V}^\B\to M_3((\mathcal V^\DR)^{\otimes2\wedge})$ 
whose equality is stated by the commutativity of (A3); its equality is then 
a consequence of the fact that the images of the generator $1_\B\in
\hat{\mathcal M}^\B$ coincide, which is itself a consequence of the definition of 
$R_{(\mu,\Phi)}$. 

The commutativity of (M4) is a consequence of two equalities, one in 
$M_{1\times 3}((\mathcal V^\DR)^{\otimes 2\wedge})$ and the other in 
$M_{3\times 1}(\mathcal M^\DR[e_1^{-1}]^{\otimes 2\wedge})$. 
The first equality is a part of the proof of the commutativity of (A4). 
The second equality follows from explicit computation based on \eqref{eq:R:mu:Phi}
and on the hexagon identities satisfied by $\Phi$. One easily 
derives the commutativity of the announced diagram.
\hfill\qed\medskip 

\subsection{Inclusion of the scheme of associators in the double shuffle scheme} 

\begin{thm}
Let $\mu\in\mathbf k^\times$. Then $M_\mu(\mathbf k)\subset
\mathsf{DMR}_\mu(\mathbf k)$. 
\end{thm}

\proof Let $\Phi\in M_\mu(\mathbf k)$. One has $(\Phi|e_0)=(\Phi|e_1)=0$ and $(\Phi|e_0e_1)=\mu^2/24$. 
Applying the equality $\hat\Delta^{\mathcal M,\DR}\circ {}^\Gamma\!\mathrm{comp}_{(\mu,\Phi)}^{\mathcal M,(10)}
=(^\Gamma\!\mathrm{comp}_{(\mu,\Phi)}^{\mathcal M,(10)})^{\otimes2}\circ\hat\Delta^{\mathcal M,\B}$ 
from Theorem \ref{thm:part:A:2804}, 2) to $1_\B\in\hat{\mathcal M}^\B$, and using 
$\hat\Delta^{\mathcal M,\B}(1_\B)=1_\B^{\otimes2}$, one obtains the group-likeness of 
$^\Gamma\!\mathrm{comp}_{(\mu,\Phi)}^{\mathcal M,(10)}(1_\B)$ for $\hat\Delta^{\mathcal M,\DR}$. 
One computes $^\Gamma\!\mathrm{comp}_{(\mu,\Phi)}^{\mathcal M,(10)}(1_\B)=
(\Gamma_\Phi(-e_1)^{-1}\Phi)\cdot 1_\DR$, which implies the result. 
\hfill\qed\medskip

\section{Bitorsor structure on the double shuffle torsor}

\subsection{The torsors 
$_{\mathsf{DMR}_0(\mathbf k)}\mathsf{DMR}_\mu(\mathbf k)$
and $_{\mathsf{DMR}^\DR(\mathbf k)}\mathsf{DMR}^{\DR,\B}(\mathbf k)$}

\begin{defn}
A {\it torsor} $_GX$ is the data of a group $G$, of a nonempty set $X$, and of a free
and transitive action of $G$ on $X$. 
\end{defn}
The left regular action of a group $G$ on itself gives rise to the
{\it trivial torsor} $_GG$. 

\begin{defn}
A torsor $_{G'}X'$ is called a subtorsor of the torsor $_GX$ iff
$G'$ (resp. $X'$) is a subgroup (resp. subset) of $G$ (resp. of $X$)
and if the action of $G'$ on $X'$ is compatible with the action of $G$ on $X$. 
\end{defn}

\begin{thm} (\cite{Rac}, \S3.2.3)
$\mathsf{DMR}_0(\mathbf k)$ is a subgroup of $(\mathcal G(\hat{\mathcal V}^\DR)
,\circledast)$, and for any $\mu\in\mathbf k^\times$, 
$_{\mathsf{DMR}_0(\mathbf k)}\mathsf{DMR}_\mu(\mathbf k)$ is a subtorsor of 
$_{\mathcal G(\hat{\mathcal V}^\DR)}\mathcal G(\hat{\mathcal V}^\DR)$. 
\end{thm}

\begin{defn}
One sets $\mathsf{DMR}^\DR(\mathbf k):=
\mathbf k^\times\times\mathsf{DMR}_0(\mathbf k)\subset
G^\DR(\mathbf k)$, $\mathsf{DMR}^{\DR,\B}(\mathbf k)
:=\{(\mu,g)|\mu\in\mathbf k^\times,g\in\mathsf{DMR}_\mu(\mathbf k)\}
\subset G^\DR(\mathbf k)$. 
\end{defn}

\begin{lem} (see \cite{EF2}, Lemma 2.13)
$\mathsf{DMR}^\DR(\mathbf k)$ is a subgroup of 
$(G^\DR(\mathbf k),\circledast)$, and 
$_{\mathsf{DMR}^\DR(\mathbf k)}\mathsf{DMR}^{\DR,\B}(\mathbf k)$ 
is a subtorsor of $_{G^\DR(\mathbf k)}G^\DR(\mathbf k)$. 
\end{lem}

\subsection{Relation of 
$_{\mathsf{DMR}^\DR(\mathbf k)}\mathsf{DMR}^{\DR,\B}(\mathbf k)$
with a stabilizer subtorsor of $_{G^\DR(\mathbf k)}G^\DR(\mathbf k)$}

\begin{lem} (see \cite{EF2}, Lemma 2.3) If $_{G_1}X_1$ and $_{G_2}X_2$ are 
subtorsors of the torsor $_GX$ such that $X_1\cap X_2\neq\emptyset$, then 
$_{G_1\cap G_2}X_1\cap X_2$ is a subtorsor of $_GX$, called the intersection of both subtorsors. 
\end{lem}

\begin{lem} (see \cite{EF2}, Lemma 2.10) Set 
$G^\DR_{\mathrm{quad}}(\mathbf k):=\{(\mu,g)\in
G^\DR(\mathbf k)|(g|e_0)=(g|e_1)=(g|e_0e_1)=0\}$
and $G^{\DR,\B}_{\mathrm{quad}}(\mathbf k):=\{(\mu,\Phi)\in
G^\DR(\mathbf k)|(\Phi|e_0)=(\Phi|e_1)=0,(\Phi|e_0e_1)=\mu^2/24\}$, then 
$_{G^\DR_{\mathrm{quad}}(\mathbf k)}G^{\DR,\B}_{\mathrm{quad}}(\mathbf k)$
is a subtorsor of $_{G^\DR_{\mathrm{quad}}(\mathbf k)}
G^\DR_{\mathrm{quad}}(\mathbf k)$. 
\end{lem}

\begin{defn}
An action of a torsor $_GX$ on a pair of isomorphic $\mathbf k$-modules
$(V,V')$ is the data of a group morphism 
$\rho:G\to\mathrm{Aut}_{\mathbf k\text{-mod}}(V)$ and of a map 
$\rho':X\to\mathrm{Iso}_{\mathbf k\text{-mod}}(V',V)$, such that 
$\rho'(g\cdot x)=\rho(g)\circ\rho'(x)$ for $g\in G$, $x\in X$. 
\end{defn}

One proves: 
\begin{lem} (see \cite{EF2} Lemma 2.6)
Let $_GX$ be a torsor and $(\rho,\rho')$ be an action on the pair 
$(V,V')$ of isomorphic $\mathbf k$-modules. Let $(v,v')\in V\times V'$. 
Let $\mathrm{Stab}(v):=\{g\in G|\rho(g)(v)=v\}$ and $\mathrm{Iso}(v',v):=
\{x\in X|\rho'(x)(v')=v\}$. If $\mathrm{Iso}(v',v)$ is nonempty, then 
$_{\mathrm{Stab}(v)}\mathrm{Iso}(v',v)$ is a subtorsor of $_GX$. 
\end{lem}

Set $_GX:=_{G^\DR(\mathbf k)}G^\DR(\mathbf k)$ and $V^\omega:=
\mathrm{Hom}_{\mathbf k\text{-mod}_{\mathrm{top}}}(\hat{\mathcal M}^\omega,
(\mathcal M^\omega)^{\otimes2\wedge})$ for $\omega\in\{\DR,\B\}$; here 
$\mathbf k\text{-mod}_{\mathrm{top}}$ is the category of topological 
$\mathbf k$-modules, i.e., $\mathbf k$-modules equipped with a decreasing 
$\mathbb Z_{\geq 0}$-filtration, separated and complete for the corresponding topology. 
An action of $_GX$ on $(V^\DR,V^\B)$ is given by $\rho : (\mu,g)\mapsto 
(V^\DR\ni f^\DR\mapsto ( ^\Gamma\!\mathrm{aut}_{(\mu,g)}^{\mathcal M,\DR,(10)}
)^{\otimes2}\circ f^\DR\circ ( ^\Gamma\!\mathrm{aut}_{(\mu,g)}^{\mathcal M,\DR,(10)}
)^{-1}\in V^\DR)$ and $\rho' : (\mu,g)\mapsto (V^\B\ni f^\B\mapsto 
( ^\Gamma\!\mathrm{comp}_{(\mu,g)}^{\mathcal M,(10)})^{\otimes2}\circ 
f^\B\circ ( ^\Gamma\!\mathrm{comp}_{(\mu,g)}^{\mathcal M,(10)})^{-1}\in V^\DR)$. 

The stabilizer subtorsor relative to the pair of vectors
$(\hat\Delta^{\mathcal M,\DR},\hat\Delta^{\mathcal M,\B})$
is denoted $_{\mathrm{Aut}(\hat\Delta^{\DR})}\mathrm{Iso}
(\hat\Delta^{\B/\DR})$. 

\begin{thm} (see \cite{EF2}, Theorem 3.1)
The subtorsor $_{\mathsf{DMR}^\DR(\mathbf k)}
\mathsf{DMR}^{\DR,\B}(\mathbf k)$ of 
$_{G^\DR(\mathbf k)}G^\DR(\mathbf k)$
 coincides with the intersection of the subtorsors
$_{\mathrm{Aut}(\hat\Delta^{\mathcal M,\DR})}\mathrm{Iso}
(\hat\Delta^{\mathcal M,\B/\DR})$ and $  _{G^\DR_{\mathrm{quad}}(\mathbf k)}
G^{\DR,\B}_{\mathrm{quad}}(\mathbf k)$. 
\end{thm}

{\it Sketch of proof.} It follows from the proof of Theorem 3.3 that $(\mu,\Phi)\in\mathrm{Iso}
(\hat\Delta^{\mathcal M,\B/\DR})$ implies $(\Gamma_\Phi(-e_1)^{-1}\Phi)\cdot 1_\DR
\in\mathcal G(\hat{\mathcal M}^\DR)$. Therefore $\mathrm{Iso}
(\hat\Delta^{\mathcal M,\B/\DR})\cap G_{\mathrm{quad}}^{\DR,\B}(\mathbf k)\subset 
\mathsf{DMR}^{\B,\DR}(\mathbf k)$. Both sides of this inclusion are subtorsors
of $_{G^\DR(\mathbf k)}G^\DR(\mathbf k)$, with underlying groups 
$\mathrm{Aut}(\hat\Delta^{\mathcal M,\DR})\cap G^\DR_{\mathrm{quad}}(\mathbf k)$ 
and $\mathsf{DMR}^{\DR}(\mathbf k)$. It follows from \cite{EF0} that these 
subgroups of $G^\DR(\mathbf k)$ are equal, which implies the equality of both torsors.   \hfill\qed\medskip

\subsection{Computation of $\mathrm{Aut}_{\mathsf{DMR}^\DR(\mathbf k)}(
\mathsf{DMR}^{\DR,\B}(\mathbf k))$ and $\mathrm{Aut}_{\mathsf{DMR}_0(
\mathbf k)}(\mathsf{DMR}_\mu(\mathbf k))$}

\subsubsection{Group corresponding to a torsor}

For $_GX$ a torsor, let $\mathrm{Aut}_G(X)$ be the group of 
right-acting permutations of $X$ which commute with the action 
of $G$. This group acts simply and transitively 
on $X$. We will call it the {\it group corresponding to} $_GX$. 

Note that the choice of an element of $X$ induces an isomorphism 
between $G$ and $\mathrm{Aut}_G(X)$, which however gets composed with an 
inner automorphism upon change of the element.

\begin{lem}\label{lemma:aut:0805}
(a) If $_{G'}X'$ is a subtorsor of $_GX$, then $\mathrm{Aut}_{G'}(X')$
is canonically a subgroup of $\mathrm{Aut}_G(X)$. 

(b) If $_{G_1}X_1$ and $_{G_2}X_2$ are subtorsors of $_GX$ with 
$X_1\cap X_2\neq\emptyset$, then $\mathrm{Aut}_{G_1\cap G_2}(X_1\cap X_2)=
\mathrm{Aut}_{G_1}(X_1)\cap\mathrm{Aut}_{G_2}(X_2)$ (equality of subgroups of 
$\mathrm{Aut}_G(X)$). 
\end{lem} 

\proof For $E$ a set, denote by $S_E$ the group of permutations of $E$, and if $e\in E$, let 
$S_{E,e}$ be the subgroup of permutations which take $e$ to itself. Then $X'$ may be 
viewed as an element of the quotient $G'\setminus X$, and there is a natural group morphism 
$\mathrm{Aut}_G(X)\to S_{G'\setminus X}$. Let $\mathrm{Aut}_G(X,X')$ be the preimage of 
$S_{G'\setminus X,X'}$. The natural map $\mathrm{Aut}_G(X,X')\to\mathrm{Aut}_{G'}(X')$
is a group isomorphism. The result follows from the diagram $\mathrm{Aut}_{G'}(X')
\simeq\mathrm{Aut}_G(X,X')\subset\mathrm{Aut}_G(X)$. This implies (a). (b) then follows from 
$\mathrm{Aut}_G(X,X_1\cap X_2)=\mathrm{Aut}_G(X,X_1)\cap\mathrm{Aut}_G(X,X_2)$. 
\hfill\qed\medskip 

\subsubsection{The group $G^\B(\mathbf k)$ and its actions}

Let $\mathcal G(\hat{\mathcal V}^\B)$ be the set of group-like of 
elements of $\hat{\mathcal V}^\B$ for $\hat\Delta^{\mathcal V,\B}$. 
Let $G^\B(\mathbf k):=\mathbf k^\times\times\mathcal G(\hat{\mathcal V}^\B)$. 
For $(\mu,g)\in G^\B(\mathbf k)$, let $\mathrm{aut}^{\mathcal V,\B,
(1)}_{(\mu,g)}$ be the automorphism of the topological $\mathbf k$-algebra 
$\hat{\mathcal V}^\B$ given by $X_0\mapsto gX_0^\mu g^{-1}$, $X_1\mapsto
X_1^\mu$, where $a\mapsto a^\mu$ is the self-map of $\mathcal G(
\hat{\mathcal V}^\B)$ given by $a^\mu:=\mathrm{exp}(\mu\mathrm{log}(a))$. 
Let $\mathrm{aut}_{(\mu,g)}^{\mathcal V,\B,(10)}$ be the topological 
$\mathbf k$-module automorphism 
of $\hat{\mathcal V}^\B$ defined by $\mathrm{aut}_{(\mu,g)}^{\mathcal V,
\B,(10)}(a):=
\mathrm{aut}_{(\mu,g)}^{\mathcal V,\B,(1)}(a)\cdot g$ for any $a\in
\hat{\mathcal V}^\B$. As in the de Rham situation, 
there is a unique 
topological $\mathbf k$-module automorphism $\mathrm{aut}_{(\mu,g)}^{
\mathcal M,\B,(10)}$ of 
$\hat{\mathcal M}^\B$, such that $\mathrm{aut}_{(\mu,g)}^{\mathcal M,
\B,(10)}(a\cdot 1_\B)
=\mathrm{aut}_{(\mu,g)}^{\mathcal V,\B,(10)}(a)\cdot 1_\B$ for any 
$a\in\hat{\mathcal V}^\B$.

One checks that $(\mu,g)\circledast(\mu',g'):=(\mu\mu',\mathrm{aut}_{(\mu,g)}^{\mathcal V,\B,(10)}(g'))$
equips $G^\B(\mathbf k)$ with a group structure, of which $\mathcal G(\hat{\mathcal V}^\B)$ is a subgroup. 

Then the map taking $(\mu,g)$ to $\mathrm{aut}_{(\mu,g)}^{\mathcal V,\B,(10)}$ 
(resp. $\mathrm{aut}_{(\mu,g)}^{\mathcal W,\B,(1)}$, $\mathrm{aut}_{(\mu,g)}^{\mathcal M,\B,(10)}$) is a 
group morphism from $(G^\B(\mathbf k),\circledast)$ to 
$\mathrm{Aut}_{\mathbf k\text{-alg}_{\mathrm{top}}}(\hat{\mathcal V}^\B)$
(resp. $\mathrm{Aut}_{\mathbf k\text{-mod}_{\mathrm{top}}}(\hat{\mathcal V}^\B)$,
$\mathrm{Aut}_{\mathbf k\text{-mod}_{\mathrm{top}}}(\hat{\mathcal M}^\B)$). 

For $g\in\hat{\mathcal V}^\B$, let $\Gamma_g(t):=\mathrm{exp}(\sum_{n\geq1}(-1)^{n+1}(g|
(\mathrm{log}X_0)^{n-1}\mathrm{log}X_1)t^n/n)\in\mathbf k[[t]]^\times$, where $w\mapsto (g|w)$ is the map 
$\{\mathrm{log}X_0,\mathrm{log}X_1\}^*\to\mathbf k$ such that $g=\sum_{w\in
\{\mathrm{log}X_0,\mathrm{log}X_1\}^*}(g|w)w$. 

As in the de Rham case, the map $\Gamma:G^\B(\mathbf k)\to(\hat{\mathcal W}^\B)^\times$, 
$(\mu,g)\mapsto \Gamma_g^{-1}(-\mathrm{log}X_1)$ satisfies a cocycle identity, which implies that the map taking 
$(\mu,g)$ to ${}^\Gamma\!\mathrm{aut}_{(\mu,g)}^{\mathcal M,\B,(10)}:=\mathrm{Ad}_{\Gamma(\mu,g)}\circ
\mathrm{aut}_{(\mu,g)}^{\mathcal M,\B,(10)}$ is a group morphism. 

There is a unique isomorphism $i^{\mathcal V}:\hat{\mathcal V}^\B\to
\hat{\mathcal V}^\DR$, induced by $X_i\mapsto\mathrm{exp}(e_i)$ for $i=0,1$. 
It induces a group isomorphism $i^G:(G^\B(\mathbf k),\circledast)\to
(G^\DR(\mathbf k),\circledast)$. 

\subsubsection{Subgroups of $G^\B(\mathbf k)$}

One checks that $G^\B_{\mathrm{quad}}(\mathbf k):=\{(\mu,g)\in G^\B(\mathbf k)|\mu^2=
1+24(g|\mathrm{log}X_0\mathrm{log}X_1)\}$ is a subgroup of $(G^\B(\mathbf k),\circledast)$ 
(see \cite{EF3}, Lemma 3.2). On the other hand, its follows from the group morphism property of 
$(\mu,g)\mapsto ^\Gamma\!\mathrm{aut}_{(\mu,g)}^{\mathcal M,\B,(10)}$ 
that $\mathrm{Aut}(\hat\Delta^{\mathcal M,\B}):=\{(\mu,g)\in G^\B(\mathbf k)|
\hat\Delta^{\mathcal M,\B}\circ\  ^\Gamma\!\mathrm{aut}_{(\mu,g)}^{\mathcal M,\B,(10)} 
=( ^\Gamma\!\mathrm{aut}_{(\mu,g)}^{\mathcal M,\B,(10)})^{\otimes2}
\circ\hat\Delta^{\mathcal M,\B}\}$ is a also a subgroup of $(G^\B(\mathbf k),\circledast)$. 

We then define $\mathsf{DMR}^\B(\mathbf k):=\mathrm{Aut}(\hat\Delta^{\mathcal M,\B})\cap 
G^\B_{\mathrm{quad}}(\mathbf k)$ and $\mathsf{DMR}^\B_0(\mathbf k)$ as the intersection of 
$\mathsf{DMR}^\B(\mathbf k)$ with $\mathcal G(\hat{\mathcal V}^\B) \subset G^\B(\mathbf k)$. 
These are subgroups of $(G^\B(\mathbf k),\circledast)$. 

One proves:
\begin{prop} (\cite{EF3})
One has $\mathsf{DMR}^\B(\mathbf k)=\{(\mu,g)\in G^\B(\mathbf k)|
(g|\mathrm{log}X_0)=(g|\mathrm{log}X_1)=0, \mu^2=1+24(g|\mathrm{log}X_0\mathrm{log}X_1), 
(\Gamma_{g}(-\mathrm{log}X_1)^{-1}\cdot g)\cdot 1_\B\in\mathcal G(\hat{\mathcal M}^\B)\}$
and $\mathsf{DMR}^\B_0(\mathbf k)=\{g\in\mathcal G(\hat{\mathcal V}^\B)|(1,g)\in 
\mathsf{DMR}^\B(\mathbf k)\}$. 
\end{prop} 

\subsubsection{Computation of groups corresponding to torsors}

\begin{thm}\label{thm:comp:aut:0805}
There are compatible group isomorphisms of 
$\mathrm{Aut}_{\mathsf{DMR}^\DR(\mathbf k)}(\mathsf{DMR}^{\DR,\B}(\mathbf k))$ 
with $\mathsf{DMR}^\B(\mathbf k)$ and for any $\mu\in\mathbf k^\times$, of 
$\mathrm{Aut}_{\mathsf{DMR}_0(k)}(\mathsf{DMR}_\mu(\mathbf k))$ with 
$\mathsf{DMR}^\B_0(\mathbf k)$. 
\end{thm}

\proof 
The group $\mathrm{Aut}_{G^\DR(\mathbf k)}(G^\DR(\mathbf k))$ corresponding to the trivial 
torsor $_{G^\DR(\mathbf k)}G^\DR(\mathbf k)$ is equal to $G^\DR(\mathbf k)$. 
According to Lemma \ref{lemma:aut:0805} (a), $\mathrm{Aut}_{G^\DR_{\mathrm{quad}}(\mathbf k)}
(G^{\DR,\B}_{\mathrm{quad}}(\mathbf k))$
is then a subgroup of $G^\DR(\mathbf k)$. One checks that its image under the isomorphism $(i^G)^{-1}$
is the subgroup $G^\B_{\mathrm{quad}}(k)$. In the same way, $\mathrm{Aut}_{\mathrm{Aut}
(\hat\Delta^{\mathcal M,\DR})}\mathrm{Iso}(\hat\Delta^{\mathcal M,\B/\DR})$  
is a subgroup of $G^\DR(\mathbf k)$, whose image under $(i^G)^{-1}$
is the subgroup $\mathrm{Aut}(\hat\Delta^{\mathcal M,\B})$. 
It then follows from Lemma \ref{lemma:aut:0805}, (b) that the image under $(i^G)^{-1}$
of $\mathrm{Aut}_{\mathsf{DMR}^\DR(\mathbf k)}(\mathsf{DMR}^{\DR,\B}(\mathbf k))$
is equal to $\mathsf{DMR}^\B(\mathbf k)$. This implies the first statement. 
The second statement follows from the fact that the natural map 
$\mathsf{DMR}^{\DR,\B}(\mathbf k)\to\mathbf k^\times$ is compatible 
with the group morphisms $\mathsf{DMR}^\omega(\mathbf k)\to
\mathbf k^\times$, $\omega\in\{\B,\DR\}$, and with the left and right actions. 
\hfill\qed\medskip 

\subsection{Bitorsor structures}

A {\it bitorsor} $_GX_H$ is a triple $(G,X,H)$ such that $_GX$ is a torsor, 
and $H$ is a group equipped with a simple and transitive right action on $X$, 
commuting with that of $G$. Like torsors, bitorsors from a category, and 
a category equivalence from torsors to bitorsors is given by 
$_GX\mapsto _GX_{\mathrm{Aut}_G(X)}$. 
Theorem \ref{thm:comp:aut:0805} may therefore be interpreted as an 
explicitation of the bitorsors corresponding to the torsors
$_{\mathsf{DMR}^\DR(\mathbf k)}\mathsf{DMR}^{\DR,\B}(\mathbf k)$ 
and $_{\mathsf{DMR}_0(\mathbf k)}\mathsf{DMR}_\mu(\mathbf k)$.

\end{document}